\documentclass[reqno, 11pt]{amsart}
\usepackage{mathrsfs}
\usepackage{amsmath,amsthm,amssymb,amscd}
\usepackage{latexsym}
\usepackage{multirow} 
\usepackage{color}
\usepackage[colorlinks=true]{hyperref}
\hypersetup{urlcolor=blue, citecolor=blue,linkcolor=blue}

\newtheorem{theorem}{Theorem}[section]

\newtheorem{lemma}[theorem]{Lemma}

\newtheorem{example}[theorem]{Example}

\theoremstyle{definition} \theoremstyle{remark}
\numberwithin{equation}{section}

\def\C{\mathbb C}
\def\R{\mathbb R}

\def\Z{\mathbb Z}

\def\Z{\mathbb Z}

\def\P{\mathbf{P}}

\def\<{\langle}
\def\>{\rangle}

\def\diag {\mathrm{diag}}

\newcommand{\tr}{\operatorname{tr}}
\newcommand{\rank}{\operatorname{rank}}

\def\F{\mathbb{F}}
\renewcommand{\i}{\mathbf{i}}
\newcommand{\be}[2]{\begin{#1} #2 \end{#1}}  
\newcommand{\mtx}[1]{\begin{pmatrix} #1 \end{pmatrix}} 
\newcommand{\wt}[1]{\widetilde{#1}}
\newcommand{\wh}[1]{\widehat{#1}}
\newcommand{\Ima}{\operatorname{Im}}

\newcommand{\RS}{\mathbf{RS}}
\newcommand{\CS}{\mathbf{CS}}
\newcommand{\DS}{\mathbf{DS}}
\newcommand{\spec}{\operatorname{spec}}

\begin{document}
	
	\title{Multiplicative trace and spectrum preservers  on  stochastic matrices}

	\author[M. Tsai] {Ming-Cheng Tsai}
	\address{General Education Center,
		Taipei University of Technology, Taipei 10608, Taiwan}
	\email{mctsai2@mail.ntut.edu.tw}
	
	\author[H. Huang] {Huajun Huang}
	\address{Department of Mathematics and Statistics, Auburn University,
		AL 36849--5310, USA}
	\email{huanghu@auburn.edu}

	

	\begin{abstract}
		We characterize maps $\phi_i: \mathcal{S} \to \mathcal{S}$, $i=1, \ldots, m$ and $m\ge 1$, that have the multiplicative spectrum or trace preserving property: 
		\begin{eqnarray*}
			\spec (\phi_1(A_1)\cdots \phi_m(A_m)) &=& \spec (A_1\cdots A_m),\quad\text{or}
			\\
			\tr (\phi_1(A_1)\cdots \phi_m(A_m)) &=& \tr (A_1\cdots A_m), 
		\end{eqnarray*}
		where  $\mathcal{S}$ is the set  of $n\times n$ doubly stochastic, row stochastic, or column stochastic matrices, or the space spanned by one of  these sets. 
		Linearity is assumed when $m=1$. We show that every stochastic matrix contains a real doubly stochastic component that carries the spectral information. 
		In consequence, the multiplicative spectrum or trace preservers on these sets $ \mathcal{S} $ are   linked to the 
		corresponding preservers on the space of doubly stochastic matrices. 
		Moreover, when $m\ge 3$,  multiplicative trace preservers always coincide with multiplicative spectrum preservers. 
	\end{abstract}

	\maketitle
	
	\section{Introduction}

	Let $M_{m,n}(\F)$ be the set of $m\times n$ matrices over a field $\F$, and write $M_{n}(\F)=M_{n,n}(\F)$.
	Denote by $\P_n$, $\DS_n$, $\RS_n$, and $\CS_n$ the sets of $n\times n$ permutation, doubly stochastic, row stochastic, and column stochastic matrices, respectively.
	For a subset $\mathcal{S}$ of a vector space over $\F$, let $\<\mathcal{S}\>_{\F}$ denote the subspace spanned by $\mathcal{S}$, and write $\<\mathcal{S}\>=\<\mathcal{S}\>_{\R}$. For a matrix $A$, let $A^t$ denote its transpose.

	The main goal of this paper is to determine the maps $\phi_i: \mathcal{S}\to \mathcal{S}$ ($i=1,\ldots,m$), where $\mathcal{S}\in\{\DS_n, \RS_n, \CS_n, \<\DS_n\>, \<\RS_n\>, \<\CS_n\>\}$,
	that  have    either {\em the multiplicative trace preserving property}:
	\begin{equation}\label{trace preserving eq}
		\tr (\phi_1(A_1)\cdots \phi_m(A_m))=\tr (A_1\cdots A_m)
	\end{equation}
	or   {\em the multiplicative  spectrum preserving property}:
	\begin{equation}\label{spectrum preserving eq}
		\spec (\phi_1(A_1)\cdots \phi_m(A_m))=\spec (A_1\cdots A_m)
	\end{equation}
	for all  $A_1,\ldots,A_m\in \mathcal{S}$. Here $\spec (A)$ denotes the spectrum of $A$ counting multiplicities, and linearity of   $\phi_1$ is assumed  when $m=1$.   A list of maps $\phi_1,\ldots,\phi_m $ is called a  {\em multiplicative trace} (resp.  {\em spectrum}) {\em preserver}  on $\mathcal{S}$ 
	if it satisfies equation \eqref{trace preserving eq} (resp. \eqref{spectrum preserving eq}). It is obvious that   \eqref{spectrum preserving eq} implies \eqref{trace preserving eq}.
	Interestingly, we show that \eqref{trace preserving eq} also implies \eqref{spectrum preserving eq}, except for $m\in\{1,2\}$ and $\mathcal{S}\in{\<\DS_n\>, \<\RS_n\>, \<\CS_n\>}$.

	Preserver problems study maps between spaces or sets that preserve specified subsets  or relations. In particular, linear preserver problems form a central and active area of research in matrix theory (e.g., \cite{Li01, Li92, Molnar, Pierce}). 
	Stochastic matrices arise naturally in many applied fields, such as probability theory, Markov chain analysis,  and statistical physics. The characterizations of maps that preserve trace or spectrum on stochastic matrices 
	can help determine transformations of Markov chains that preserve transition probabilities’ equilibrium properties  or 
	long-term statistical behaviors in probability and physics models.

	Trace preservers and spectrum preservers are closely related. 
	On one hand, spectrum preservers are necessarily trace preservers. On the other hand, if two $n\times n$ matrices $A$ and $B$ satisfy $\tr(A^k)=\tr(B^k)$ for all $k\in\Z^+$, then the symmetric functions of the eigenvalues of $A$ and $B$ coincide, and hence  $\spec(A)=\spec(B)$.

	Marcus and Moyls initiated the  study of spectrum preserving problem in 1959 \cite{Ma59}, showing that a linear map $\phi: M_n(\C)  \to M_n(\C) $ satisfying 
	\begin{equation}
		\label{spectrum}
		\spec(\phi(A))=\spec(A) 
	\end{equation}
	for $A\in M_n(\C)$ 
	must take the form $\phi(A)=SAS^{-1}$ or $\phi(A)=SA^tS^{-1}$ for  an invertible matrix $S\in M_n(\C)$. 
	In 1974, Minc \cite{Minc} proved that if a linear operator $\phi$ on 
	$M_n(\C) $ preserves nonnegative matrices and satisfies   \eqref{spectrum}  for all nonnegative matrices, then
	$\phi(A)=PAP^{-1}$ or $\phi(A)=PA^tP^{-1}$ where $P$ is the product of a permutation matrix and a  positive diagonal matrix. 
	Later, some researchers extended the research direction to multiplicative spectrum preservers
	(e.g., \cite{ChanLiSze, ClarkLiRodman, Hou08, Hou10, Molnar2, Molnar3}). In 2003, Molnár \cite{Molnar3} showed that if a map $\phi: M_n(\C) \to M_n(\C) $ satisfies
	\begin{equation}
		\label{spectrum5}
		\spec(\phi(A)\phi(B))=\spec(AB) 
	\end{equation}
	for $A, B \in M_n(\C)$, 
	then there exists
	an invertible matrix $S$ and $\lambda\in \{-1, 1\}$  such that $\phi(A)=\lambda SAS^{-1}$ or   $\phi(A)=\lambda SA^tS^{-1}$. In 2008, Clark, Li, and Rodman \cite{ClarkLiRodman} characterized preservers of the spectral radius, spectrum and peripheral spectrum of some products or  Jordan triple products on nonnegative matrices. More recently, several researchers have studied related problems involving two or more maps preserving certain spectra on various matrix spaces (e.g., \cite{Abdelali, Bourhim, Huang23, Mi09}). In particular, Abdelali and Aharmim \cite{Abdelali} investigate the maps $\phi: M_{m,n}(\mathbb{C})\rightarrow M_{p,q}(\mathbb{C})$ and $\psi: M_{n,m}(\mathbb{C})\rightarrow M_{q,p}(\mathbb{C})$ satisfying 
	\be{equation}
	{\label{spectrum4}
		\sigma(\phi(A)\psi(B))=\sigma(AB)
	}
	for all $A\in
	M_{m,n}(\mathbb{C}), B\in M_{n,m}(\mathbb{C}),$ where $\sigma (A)$ denotes the spectrum of $A$ without counting multiplicities.

	The known earliest classical result about   multiplicative trace preservers appeared in Wigner's   theorem {\cite[p.12]{Molnar}}. It states that if a bijective operator $\phi$  on the set of all rank one projections acting on a Hilbert space $H$  satisfies 
	\be{equation}
	{\label{trace}
		\tr(\phi(P)\phi(Q))=\tr(PQ)
	}
	for all rank one projections $P$, $Q$, then there exists a unitary  operator $U$ on $H$ such that $\phi(P)=U^{*}PU$
	or $\phi(P)=U^{*}P^{t} U$ for all  rank one projection $P$. Since then, a number of researchers have investigated maps that preserve various trace equations (see, e.g., \cite{Huang16, Huang23, Leung, Li12, Uh63}). In 2023, Huang and Tsai  characterized the maps $\phi_i : \mathcal{S} \to \mathcal{S}$ ($i = 1, \ldots, m$) satisfying \eqref{trace preserving eq}, where $\mathcal{S}$ denotes the set of $n \times n$ general, Hermitian, or symmetric matrices over $\F$ for $m \ge 3$, and the set of positive definite or diagonal matrices over $\F$ for $m \ge 2$, with $\F = \C$ or $\R$ (see \cite{Huang23}). The authors' works indicate  that \eqref{trace preserving eq} is equivalent to \eqref{spectrum preserving eq} in most of the cases.

	In this paper, we completely characterize the lists of maps $\phi_1, \ldots, \phi_m$ on the set $\mathcal{S}\in\{\DS_n, \RS_n, \CS_n, \<\DS_n\>, \<\RS_n\>, \<\CS_n\>\}$ that satisfy \eqref{trace preserving eq} or \eqref{spectrum preserving eq}. 
	
	Section 2 contains preliminary results. Let $e=(1,\dots,1)^t\in \R^n$ and define
	$$R_n=\{e v^t:\ v\in \R^n, v^t e=0\},\qquad
	C_n=\{ve^t:\ v\in \R^n, v^t e=0\}. 
	$$
	Lemma \ref{thm: M_n=DS+L+L^t} gives the direct sums of matrix subalgebras
	$$M_n(\R)= \< \DS_n \>\oplus R_n\oplus C_n,\quad
	\< \RS_n \> =  \< \DS_n \>\oplus R_n,\quad
	\< \CS_n \> =  \< \DS_n \>\oplus C_n.
	$$
	We then determine the structures of $\<\DS_n\>$, $\<\RS_n\>$, and $\<\CS_n\>$ (Theorem \ref{thm: DS RS CS description}) and
	characterize the algebra automorphisms on these spaces (Theorem \ref{thm: stochastic algebra automorphism}). 
	Theorem \ref{thm: DS RS CS} and Lemma \ref{thm: inverse of DS}  allow us to transfer spectra and traces of maps from $\<\RS_n\>$ and  $\<\CS_n\>$ to $\<\DS_n\>$. Consequently, Theorem \ref{thm: DS RS preserver relation} establishes the connection between multiplicative spectrum (resp. trace) preservers on $\<\RS_n\>$ and on $\<\DS_n\>$. Likewise for those of $\<\CS_n\>$ and $\<\DS_n\>$.  
	Lemma \ref{thm: product permutation} shows that a product of doubly (resp. row, column) stochastic matrices is a permutation matrix if and only if 
	each product component is. 
	
	Sections 3 and 4 provide the detailed characterizations of multiplicative trace preservers and multiplicative  spectrum preservers on $\mathcal{S}\in\{\DS_n, \RS_n, \CS_n, \<\DS_n\>, \<\RS_n\>, \<\CS_n\>\}$. 
	
	For $\mathcal{S}\in\{\DS_n, \RS_n, \CS_n\}$, we determine the maps $\phi_1, \dots, \phi_m$ satisfying \eqref{spectrum preserving eq} when $m=1$, and those satisfying \eqref{trace preserving eq} or \eqref{spectrum preserving eq} when $m\ge 2$.
	These results  are summarized below.
	
	\begin{itemize}
		\item   $\mathcal{S}=\DS_n$:
		\begin{enumerate}
			\item
			For $m=1$ (Theorem \ref{thm: DS_n one map spectrum}), 
			$\phi_1(A)=PAP^t$ or  $\phi_2(A)=PA^tP^t$ for certain  $P\in \P_n$. 
			\item
			For $m=2$ (Theorem \ref{thm: DS_n two maps preserving trace of product}), 
			there exist  $P, Q\in \P_n$ such that
			$\be{cases}{\phi_1(A)=PAQ^t \\ \phi_2(A)=QAP^t}$
			or
			$\be{cases}{\phi_1(A)=PA^t Q^t \\ \phi_2(A)=QA^t P^t}.$
			\item For $m\ge 3$  (Theorem \ref{thm: DS_n m maps preserving trace of product}),
			there exist permutation matrices $P_1,\ldots,P_m, P_{m+1}=P_1\in \P_n$ such that $\phi_i(A)=P_iAP_{i+1}^t$ for $i=1,\ldots, m.$
		\end{enumerate}
		Besides, Example \ref{RS_n: preserver in <RS_n> not RS_n} shows that preservers on $\<\RS_n\>$ leaving $\DS_n$ invariant admit more forms than those on $\RS_n$.
		
		\item $\mathcal{S}=\RS_n$: 
		\begin{enumerate}
			\item
			For $m=1$ (Theorem \ref{thm: RS_n one map spectrum}), 
			\begin{enumerate}
				\item
				when $n=2$,
				$\phi(A) =A_{\DS}+sA_{R}$ for certain $s\in [-1,1]$; 
				\item
				when $n\ge 3$,  $\phi(A) =P A P^t$ for some permutation matrix $P\in \P_n$.
			\end{enumerate}
			\item
			For $m\ge 2$ (Theorem \ref{thm: RS_n maps preserving trace of products}): for $A\in\RS_n$, write $A=A_{\DS}+A_{R}$ in which $A_{\DS}\in\<\DS_n\>$ and $A_{R}\in R_n$. Then there exist $P_1,\ldots,P_m, P_{m+1}=P_1\in \P_n$ and $\gamma_i:\RS_n\to R_n$  such that $\phi_i(A)=P_iA_{\DS}P_{i+1}^t+\gamma_i(A)$, where
			$\gamma_i$ must let $\phi_i(A)\in \RS_n$ for all $A\in \RS_n$.
			
			If furthermore, $\phi_k$ is linear for certain $k\in [m]$, then $\phi_k$ has the form: 
			\begin{enumerate}
				\item
				when $n=2$,
				$\phi_k(A) =P_k (A_{\DS}+s_k A_{R})P_{k+1}^t$ for certain $s_k\in [-1,1]$; 
				\item
				when $n\ge 3$,  $\phi_k(A) =P_k A P_{k+1}^t$.
			\end{enumerate}
		\end{enumerate}
		\item   $\mathcal{S}=\CS_n$:  the results  are analogous to those of $\RS_n$ case. 
	\end{itemize}

	Let $\mathcal{S}\in\{\<\DS_n\>, \<\RS_n\>, \<\CS_n\>\}$. For $m=1$, we determine the linear map $\phi_1$ on $\mathcal{S}$ satisfying \eqref{spectrum preserving eq}. For $m\ge 2$, we characterize the maps $\phi_1, \ldots, \phi_m$ on $\mathcal{S}$  satisfying 
	\eqref{trace preserving eq} or \eqref{spectrum preserving eq}, with a mild  rank condition on the multiplicative trace preservers when $m=2$.
	These results are summarized below.
	
	\begin{itemize}
		\item $\mathcal{S}=\<\DS_n\>$:
		\begin{enumerate}
			\item
			For $m=1$ (Theorem \ref{thm: <DS> one map spectrum}), 
			\begin{enumerate}
				\item
				when $n=2$,
				$\phi(A)=PAP^{-1}$ for some  $P\in\left\{I_2, \mtx{1&0\\0 &-1} \right\}$;
				\item
				when $n\ge 3$,  $\phi_1(A)=PAP^{-1}$ or $\phi_2(A)=PA^tP^{-1}$ for an invertible $P\in\<\DS_n\>$.
			\end{enumerate}
			\item
			For $m=2$ (Theorem \ref{thm: <DS> two maps spectrum trace}), 
			\begin{enumerate}
				\item
				when $n=2$,  $\phi_1(A)=CPAP^{-1}$ or  $\phi_2(A)=PAP^{-1}C^{-1}$ for $P\in\left\{I_2, \mtx{1&0\\0 &-1} \right\}$ and an invertible $C\in\<\DS_2\>$;
				\item
				when $n\ge 3$, there exist invertible matrices $P, Q\in \<\DS_n\>$ such that 
				$
				\be{cases}{\phi_1(A)=PAQ^{-1} \\  \phi_2(A)=QAP^{-1}}$
				or 
				$
				\be{cases}{\phi_1(A)=PA^t Q^{-1} \\ \phi_2(A)=QA^t P^{-1}}.
				$
			\end{enumerate}
			\item For $m\ge 3$  (Theorem \ref{thm: <DS> m ge 3 maps spectrum trace}),
			\begin{enumerate}
				\item
				when $n=2$, we have $\phi_i(A)=PAP^{-1}Q_i$, $i=1,\ldots, m,$ where $P\in\left\{I_2, \mtx{1&0\\0 &-1} \right\}$  and $Q_1, \ldots, Q_m\in\<\DS_2\>$ satisfy that $Q_1\cdots Q_m=I_2$; 
				\item
				when $n\ge 3$,  $\phi_i(A)=P_iAP_{i+1}^{-1}$, $i=1,\ldots, m$, where  $P_1,\ldots,P_m, P_{m+1}=\<\DS_n\>$ are invertible.
			\end{enumerate}
		\end{enumerate}

		\item $\mathcal{S}=\<\RS_n\>$ or $\<\CS_n\>$:  In both cases, the characterizations (Theorems \ref{thm: span RSn one map spectrum}, \ref{thm: <RS> two maps spectrum trace}, and \ref{thm: <RS> m ge 3 maps spectrum trace}) are almost the same as those for $A\in\<\DS_n\>$, except that 
		the matrices in the domains and ranges will be replaced by their projections on $\<\DS_n\>.$
		
	\end{itemize}

	\section{Preliminary}
	
	For every positive integer $m$, denote the set
	\be{equation}{
		[m]:=\{1,2,\ldots,m\}.
	}
	Let $E_{i,j}^{(n)}$ be the  $n\times n$ matrix with the single nonzero entry 1 in the entry $(i,j)$.    When $n$ is well known, we simply write $E_{i,j}=E_{i,j}^{(n)}$. 
	When $\F$ is a field and $S$ and $T$ are subsets of some $\F$-vector spaces, 
	a function $\phi:S\to T$ is called an {\em ($\F$-)linear map} if $\phi$ can be extended to an
	$\F$-linear map between vector spaces $\langle S\rangle_{\F} \to \langle T\rangle_{\F}.$

	\subsection{Two maps preserving multiplicative trace}
	
	
	When $V_1$ (resp. $V_2$) is a subspace of $M_{m,n}(\F)$ (resp. $M_{n,m}(\F)$), the trace map 
	$\mathcal{B}:V_1\times V_2\to\F$, defined by $\mathcal{B}(A,B)=\tr(AB)$, is  
	a bilinear form. 
	The trace map is said to be \emph{non-degenerate} on $V_1\times V_2$, or be \emph{non-degenerate} on $V_1$ when $m=n$ and $V_1=V_2$, if any matrix representation of
	$\mathcal{B}$ is invertible. In particular, $\mathcal{B}$ is non-degenerate only if $\dim V_1=\dim V_2$. 
	
	The following result indicates that two maps between subsets of matrix spaces that preserve multiplicative trace can be extended to linear bijections, respectively, provided that the trace map is non-degenerate.

	\begin{theorem}[\cite{Huang23}]
		\label{thm: two maps preserving trace}
		Let  $\phi: V_1\to W_1$ and $\psi: V_2\to W_2$  be two maps between subsets of matrix spaces such that:
		\begin{enumerate}
			\item $\dim \langle V_1\rangle=\dim \langle V_2\rangle\ge \max\{\dim \langle W_1\rangle, \dim \langle W_2\rangle\}$.
			\item 
			$AB$  are well-defined square matrices for $(A,B)\in (V_1\times V_2)\cup(W_1\times W_2)$.
			\item \label{tr nonsingular}
			If $A\in \langle V_1\rangle $ satisfies $\tr(AB)=0$ for all $B\in \langle V_2\rangle$, then $A=0$.
			\item $\phi$ and $\psi$ satisfy that
			\begin{equation}
				\tr (\phi(A)\psi(B))=\tr(AB),\quad   A\in V_1,\ B\in V_2.
			\end{equation}
		\end{enumerate}
		Then
		$\dim \langle V_1\rangle =\dim \langle V_2\rangle =\dim \langle W_1\rangle =\dim \langle W_2\rangle, $
		and $\phi$ and $\psi$ can be extended uniquely to linear bijections $\widetilde\phi : \langle V_1\rangle \to  \langle W_1\rangle$ and $\widetilde\psi: \langle V_2\rangle \to  \langle W_2\rangle$, respectively, such that
		\begin{equation}\label{extended map preserve trace}
			\tr (\wt{\phi}(A)\wt{\psi}(B))=\tr(AB),\quad   A\in \langle V_1\rangle,\ B\in \langle V_2\rangle.
		\end{equation}
	\end{theorem}

	Theorem \ref{thm: two maps preserving trace} will be used in this paper for
	the cases  that $V_1=V_2$ is  $\DS_n$  or $\<\DS_n\>$. 
	The condition  (3) in Theorem \ref{thm: two maps preserving trace} does not hold when
	$V_1=V_2$ is $\<\RS_n\>$ or $\<\CS_n\>$.

	\subsection{The sets and spaces of stochastic matrices}
	
	We explore some properties of the sets $\DS_n$ of doubly stochastic matrices, $\RS_n$ of row stochastic matrices,  $\CS_n$ of column stochastic matrices,
	and the corresponding real spanning spaces $\<\DS_n\>$, $\<\RS_n\>$, and $\<\CS_n\>$.

	Let $e=(1,\cdots,1)^t$ be the all-ones vector in $\R^n$. 
	The spanning of doubly (resp. row, column) stochastic matrices are as follows:
	\be{eqnarray}{
		\label{doubly stochastic matrices}
		\<\DS_n\>&=&\{A\in M_n(\R): Ae=ce,\ e^t A=ce^t,\ c\in\R\},
		\\ 
		\label{row stochastic matrices}
		\<\RS_n\>&=&\{A\in M_n(\R): Ae=ce,\ c\in\R\},
		\\
		\label{column stochastic matrices}
		\<\CS_n\>&=&\{A\in M_n(\R): e^t A=ce^t,\ c\in\R\}.
	}
	The sets $\DS_n$, $\RS_n$, and $\CS_n$ are monoids with respect to matrix product, and   $\<\DS_n\>$, $\<\RS_n\>$, and $\<\CS_n\>$ are unital algebras. 
	
	Denote the following subspaces of $M_n(\R)$:
	\be{eqnarray}{
		R_n &=& \{e v^t:\ v\in \R^n, v^t e=0\},\\
		C_n &=& \{ve^t:\ v\in \R^n, v^t e=0\}. 
	}
	In other words, $R_n$ (resp. $C_n$) is the set of matrices in $M_n(\R)$ whose entries in each row (resp. column) are summed to zero and in each column (resp. row) are identical. 
	
	Given $A=[a_{ij}]\in M_n(\R)$, we define:
	\begin{align}
		\label{s_A}
		s_{A} &= \frac{1}{n^2} \sum_{i=1}^{n}\sum_{j=1}^{n} a_{ij},
		\\\label{r(A)}
		r(A) &= (r_1,\cdots, r_n)^t,\quad \text{where each } r_i=\frac{a_{1i}+a_{2i}+\cdots+a_{ni}}{n}-s_{A},
		\\\label{c(A)}
		c(A) &= (c_1,\cdots,c_n)^t,\quad \text{where each } c_i=\frac{a_{i1}+a_{i2}+\cdots+a_{in}}{n} -s_{A}.
	\end{align}
	Let 
	\begin{equation}
		\label{eq: decomp M_n to DS}
		A_{R}=er(A)^t,\quad A_{C}=c(A) e^t,\quad A_{\DS}=A-A_{R}-A_{C}.
	\end{equation}
	Clearly $A_{R}\in R_n$ and $A_{C}\in C_n$. 
	Moreover, the sum of entries in each row or column of $A_{\DS}$ is   $n s_{A}$, so that
	$A_{\DS}\in\<\DS_n\>$. 
	Hence
	$M_n(\R)=\< \DS_n \>+ R_n+ C_n$ such that
	\begin{equation}
		\label{A decomp stochastic}
		A=A_{\DS}+A_{R}+A_{C},\qquad A\in M_n(\R).
	\end{equation}

	\begin{lemma}\label{thm: M_n=DS+L+L^t}
		The following are vector space direct sums:
		\begin{align}
			\label{M_n direct sum of stochastic matrices}
			M_n(\R) &= \< \DS_n \>\oplus R_n\oplus C_n,
			\\\label{RS_n direct sum of stochastic matrices}
			\< \RS_n \> &=  \< \DS_n \>\oplus R_n,
			\\\label{CS_n direct sum of stochastic matrices}
			\< \CS_n \> &= \< \DS_n \>\oplus C_n.
		\end{align}
	\end{lemma}

	\begin{proof}
		We have shown that
		$M_n(\R)=\< \DS_n \>+ R_n+ C_n$. Suppose
		$0=A_1+A_2+A_3$ for $A_1\in\< \DS_n \>$, $A_2\in R_n$, and $A_3\in C_n$. Then
		$A_1, A_2\in \< \RS_n \>$, so that $A_3\in  \< \RS_n \>\cap C_n=\{0\}$. Similarly, $A_2=0$.  Hence $A_1=0$. 
		Therefore, $M_n(\R)  =  \< \DS_n \>\oplus R_n\oplus C_n.$
		
		By $\< \DS_n \>\oplus R_n\subseteq \<\RS_n\>$ and $\< \RS_n \>\cap C_n=\{0\}$, we get
		\eqref{RS_n direct sum of stochastic matrices}.  
		Likewise, \eqref{CS_n direct sum of stochastic matrices} holds. 
	\end{proof}

	It is known \cite[Theorem 29]{Brauer} that  $\spec(A)=\spec(A+X)=\spec(A+Y)$ for $A\in\DS_n$, $X\in R_n$, and $Y\in C_n$. 
	This property can be extended to a product of matrices in $\<\RS_n\>$ or $\<\CS_n\>$. 
	Choose $X\in M_{n,n-1}(\R)$ such that
	\begin{equation}
		\label{DS_n-orthogonal}
		U=\mtx{\frac{1}{\sqrt{n}}e, X}
	\end{equation} 
	is a real orthogonal matrix. 
	Given $A \in M_n(\R)$, we  write $A=A_{\DS}+er(A)^t+c(A)e^t$ as in \eqref{eq: decomp M_n to DS}. Then
	\begin{equation}\label{A decomp unitary conjugacy}
		U^{-1} A U=\mtx{\frac{1}{\sqrt{n}}e^t\\ X^t} (A_{\DS}+er(A)^t+c(A)e^t)\mtx{\frac{1}{\sqrt{n}}e, X}=
		\mtx{n s_{A}  &\sqrt{n}r(A)^tX\\ \sqrt{n} X^t c(A) &X^t A_{\DS} X}.
	\end{equation} 
	It implies the following theorem,  which extends the results in \cite[Proposition 5.2]{Alwill}.
	
	\begin{theorem}\label{thm: DS RS CS description}
		Let $U$ be given in \eqref{DS_n-orthogonal}. Then
		\begin{eqnarray} \label{<DS> unitary}
			\<\DS_n\> &=& \left\{U\mtx{a&0\\0 &B}U^{-1}: a\in\R,\ B\in M_{n-1}(\R)\right\},
			\\ \label{<RS> unitary}
			\<\RS_n\> &=& \left\{U\mtx{a&r^t\\0 &B}U^{-1}: a\in\R,\ B\in M_{n-1}(\R),\ r\in\R^{n-1}\right\},
			\\  \label{<CS> unitary}
			\<\CS_n\> &=& \left\{U\mtx{a&0\\c &B}U^{-1}: a\in\R,\ B\in M_{n-1}(\R),\ c\in\R^{n-1}\right\}.
		\end{eqnarray}
		In particular, every $A\in M_n(\R)$ has $A=A_{\DS}+A_R+A_C$ as in \eqref{A decomp stochastic}, in which:
		\begin{eqnarray}\label{A_DS}
			A_{\DS} &=& U \mtx{n s_{A}  &0\\ 0 &X^t A_{\DS} X}U^{-1},
			\\ \label{A_R}
			A_R &=& U \mtx{0  &\sqrt{n}r(A)^tX\\ 0 &0} U^{-1},
			\\ \label{A_C}
			A_C &=& U \mtx{0  &0\\ \sqrt{n}X^t c(A) &0} U^{-1}.
		\end{eqnarray}
	\end{theorem}
	
	\begin{proof}
		\eqref{A decomp unitary conjugacy} immediately implies the identities \eqref{A_DS}, \eqref{A_R}, and \eqref{A_C}. 
		By \eqref{M_n direct sum of stochastic matrices} and comparison of dimensions, we get \eqref{<DS> unitary} and
		\begin{eqnarray}
			R_n &=& \left\{U\mtx{0 &r^t\\0 &0}U^{-1}:  r\in\R^{n-1}\right\},
			\\  
			C_n &=& \left\{U\mtx{0&0\\c &0}U^{-1}:  c\in\R^{n-1}\right\},
		\end{eqnarray}
		which lead to \eqref{<RS> unitary} and \eqref{<CS> unitary}.
	\end{proof}

	\subsection{Operations on stochastic matrices}
	
	The  identities \eqref{A_DS}, \eqref{A_R}, and \eqref{A_C} imply the following results.


	\begin{theorem}\label{thm: DS RS CS}
		For any positive integer $m$ and matrices $A_1,\ldots,A_m$  in $\<\RS_n\>$  (resp. $\<\CS_n\>$), 
		\begin{equation}\label{DS parts preserve product}
			(A_1\cdots A_m)_{\DS} = (A_1)_{\DS}\cdots (A_m)_{\DS}.  
		\end{equation}
		In consequence, the following identities hold: 
		\begin{eqnarray} 
			\label{DS parts preserve spectrum}
			\spec(A_1\cdots A_m) &=& \spec((A_1)_{\DS}\cdots (A_m)_{\DS})
			\\ \label{mixed parts preserve spectrum}
			&=& \spec(A_1\cdots (A_i)_{\DS} \cdots A_m),\quad i\in [m],
			\\ \label{DS parts preserve trace}
			\tr(A_1\cdots A_m) &=& \tr((A_1)_{\DS}\cdots (A_m)_{\DS})
			\\ \label{mixed parts preserve trace}
			&=& \tr(A_1\cdots (A_i)_{\DS} \cdots A_m),\quad i\in [m]. 
		\end{eqnarray}
	\end{theorem}
	
	\begin{proof}
		Assume that $A_1,\ldots,A_m\in \<\RS_n\>$ (similarly for  the case $A_1,\ldots,A_m\in \<\CS_n\>$).
		By \eqref{A_DS} and \eqref{A_R}, every $A\in  \<\RS_n\>$ and the corresponding $A_{\DS}$ have the expressions: 
		\begin{eqnarray*}
			A =   U \mtx{n s_{A}  &\sqrt{n}r(A)^t X \\ 0 &X^t A_{\DS} X}U^{-1},
			\qquad
			A_{\DS} =  U \mtx{n s_{A}  &0\\ 0 &X^t A_{\DS} X}U^{-1}.
		\end{eqnarray*}
		We immediately get \eqref{DS parts preserve product} and the other inequalities.
	\end{proof}

	\begin {lemma} \label{thm: inverse of DS}
	Let $A\in \<\RS_n\>$ (resp. $A\in \<\CS_n\>$).  
	If $A$ is invertible, then so is $A_{\DS}$, and $A^{-1}\in \<\RS_n\>$  (resp. $A^{-1}\in \<\CS_n\>$) such that
	\begin{equation}
		\label{DS part inverse}
		(A^{-1})_{\DS}=(A_{\DS})^{-1}.
	\end{equation}
	In particular,
	if $A\in \<\DS_n\>$ is invertible, then $A^{-1}\in \<\DS_n\>$. 
\end{lemma}

\begin{proof}
	When $A\in \<\RS_n\>$ (resp. $A\in \<\CS_n\>$) is invertible, \eqref{<RS> unitary} shows 
	that $A^{-1}\in \<\RS_n\>$  (resp. $A^{-1}\in \<\CS_n\>$), and \eqref{DS parts preserve product} shows that
	$$
	AA^{-1}=I_n\quad\Longrightarrow\quad A_{\DS}(A^{-1})_{\DS}=(I_n)_{\DS}=I_n.
	$$
	We get \eqref{DS part inverse}. If, in addition, $A\in\<\DS_n\>$, then \eqref{DS part inverse} becomes
	$(A^{-1})_{\DS}=A^{-1}$, so that $A^{-1}\in\<\DS_n\>$. 
	%
\end{proof}

Theorem \ref{thm: DS RS CS} and Lemma \ref{thm: inverse of DS} provide a way to transform the spectra and traces of maps in $\<\RS_n\>$ or $\<\CS_n\>$ to those in $\<\DS_n\>$. In consequence, the relation of multiplicative spectrum (resp. trace) preservers  on $\<\RS_n\>$ and $\<\DS_n\>$ are  as follows.  
Likewise for the relation of those preservers on $\<\CS_n\>$ and $\<\DS_n\>$.

\begin{theorem}\label{thm: DS RS preserver relation}
	Let  $m$ and $n$ be  positive integers. 
	\begin{enumerate}
		\item If $\phi_i: \<\DS_n\> \to \<\DS_n\>$ ($i=1,\ldots,m$) is a  multiplicative spectrum (resp. trace) preserver on  
		$\<\DS_n\>$, and $\gamma_i: \<\RS_n\> \to R_n$ ($i=1,\ldots,m$) are arbitrary maps, then $\wt{\phi}_i: \<\RS_n\> \to \<\RS_n\>$ ($i=1,\ldots,m$),  defined by 
		\begin{equation*}
			\wt{\phi}_i (A)= \phi_i(A_{\DS})+\gamma_i(A),\qquad A\in \<\RS_n\>, 
		\end{equation*}
		is a  multiplicative spectrum (resp. trace) preserver on     $\<\RS_n\>$. Moreover, $\wt{\phi}_i$ is an extension of $\phi_i$
		from $\<\DS_n\>$ to $\<\RS_n\>$   for $i=1,\ldots,m$ if and only if each $\gamma_i$ sends $\<\DS_n\>$ to $0$. 
		\item 
		Conversely, if $\wt{\psi}_i: \<\RS_n\> \to \<\RS_n\>$ ($i=1,\ldots,m$) is a  multiplicative spectrum (resp. trace) preserver on     $\<\RS_n\>$, then $\psi_i: \<\DS_n\> \to \<\DS_n\>$ ($i=1,\ldots,m$), defined by 
		\begin{equation*}
			\psi_i (A)=\wt{\psi}_i(A)_{\DS},\qquad A\in \<\DS_n\>,
		\end{equation*}
		is a  multiplicative spectrum (resp. trace) preserver on     $\<\DS_n\>$. Moreover, when $m\ge 2$, 
		each $\wt{\psi}_i$ can be written as
		\begin{equation*}
			\wt{\psi}_i(A)= \psi_i(A_{\DS}) + \gamma_i(A) 
		\end{equation*}
		for a map $\gamma_i: \<\RS_n\> \to R_n$ defined by  $\gamma_i(A)=\wt{\psi}_i(A)_R$.  
	\end{enumerate}
\end{theorem} 

\begin{proof}
	Theorem \ref{thm: DS RS CS}   implies all statements in (1) and all but the last statement in (2). 
	
	It remains to prove that when $m\ge 2$, every  multiplicative spectrum (resp. trace) preserver on     $\<\RS_n\>$,
	$\wt{\psi}_i: \<\RS_n\> \to \<\RS_n\>$ ($i=1,\ldots,m$), can be expressed as
	$$\wt{\psi}_i(A)= \psi_i(A_{\DS}) + \gamma_i(A)=
	\wt{\psi}_i(A_{\DS})_{\DS}+\wt{\psi}_i(A)_R,$$ 
	that is, $\wt{\psi}_i(A)_{\DS} =\wt{\psi}_i(A_{\DS})_{\DS}$ for $A\in\<\RS_n\>.$ Without loss of generality, we assume that $i=1$. 
	For any $A \in\<\RS_n\>$ and $A_2\in\<\DS_n\>$, by assumption and Theorem \ref{thm: DS RS CS} we have 
	\begin{align*}
		&\tr(\wt{\psi}_1(A)_{\DS}\wt{\psi}_2(A_2)_{\DS}\wt{\psi}_3(I)_{\DS}\cdots \wt{\psi}_m(I)_{\DS})
		\\
		=&  \tr( \wt{\psi}_1(A) \wt{\psi}_2(A_2)\wt{\psi}_3(I)\cdots \wt{\psi}_m(I))
		\\
		=&  \tr(AA_2) = \tr(A_{\DS}   A_2)
		\\
		=& \tr( \wt{\psi}_1(A_{\DS})  \wt{\psi}_2(A_2)\wt{\psi}_3(I)\cdots \wt{\psi}_m(I))
		\\
		=&  \tr( \wt{\psi}_1(A_{\DS})_{\DS} \wt{\psi}_2(A_2)_{\DS}\wt{\psi}_3(I)_{\DS}\cdots \wt{\psi}_m(I)_{\DS}).
	\end{align*}
	It implies the following two things:
	\begin{enumerate}
		\item Define $\phi_2:\<\DS_n\>\to \<\DS_n\>$ by 
		$$\phi_2(A_2)=\wt{\psi}_2(A_2 )_{\DS}\wt{\psi}_3(I)_{\DS}\cdots \wt{\psi}_m(I)_{\DS}.$$
		Then for $A, A_2\in\<\DS_n\>$ we get 
		$$\tr(\psi_1(A)\phi_2(A_2)) = \tr(\wt{\psi}_i(A)_{\DS} \wt{\psi}_2(A_2 )_{\DS}\wt{\psi}_3(I)_{\DS}\cdots \wt{\psi}_m(I)_{\DS})=\tr(AA_2).$$
		Theorem \ref{thm: DS RS CS description} implies that the trace map   is non-degenerate on $\<\DS_n\>$. Theorem \ref{thm: two maps preserving trace} shows
		that both $\psi_1$ and $\phi_2$ are linear bijections on $\<\DS_n\>$. 
		\item Now for any given $A\in\<\RS_n\>$ and any $A_2\in \<\DS_n\>$ we have 
		\begin{align*}
			0 &=\tr \left ( (\wt{\psi}_1(A)_{\DS}- \wt{\psi}_1(A_{\DS})_{\DS}) \wt{\psi}_2(A_2)_{\DS}\wt{\psi}_3(I)_{\DS}\cdots \wt{\psi}_m(I)_{\DS} \right ) 
			\\ &= \tr \left ( (\wt{\psi}_1(A)_{\DS}- \wt{\psi}_1(A_{\DS})_{\DS}) \phi_2(A_2)\right). 
		\end{align*}
		Since the trace map   is non-degenerate on $\<\DS_n\>$ and $\Ima\phi_2= \<\DS_n\>$, we get $\wt{\psi}_1(A)_{\DS}=\wt{\psi}_1(A_{\DS})_{\DS}$ as desired. 
		\qedhere
	\end{enumerate} 
\end{proof}

Viewing $\DS_n$, $\RS_n$, and $\CS_n$ as monoids,
M. Alwill, C-K Li, C. Maher, and N-S Sze described  the  monoid automorphisms of $\DS_n$, $\RS_n$, and $\CS_n$ in
\cite[Theorem 6.1]{Alwill}, and the  monoid automorphisms of the sets of ``stochastic matrices'' that drop the nonnegative requirement on their entries, namely, the sets
$\{(c+1)A-c (ee^t)\mid c\in\R, A\in {\mathcal S}\}$ for $ {\mathcal S} \in \{\DS_n,\RS_n,\CS_n\}$, in \cite[Theorem 5.1]{Alwill}.  
We characterize the  algebra  automorphisms  of $ \<\DS_n\>$, $\<\RS_n\>$, and $\<\CS_n\>$ here. 


\begin{theorem}\label{thm: stochastic algebra automorphism}
	When $n\ge 2$, every  algebra  automorphism $\phi$ of ${\mathcal A}=\<\DS_n\>$, $\<\RS_n\>$, or $\<\CS_n\>$,  has the form
	\be{equation}{\label{DS RS CS algebra automorphism}
		\phi(A)=PAP^{-1},\quad A\in {\mathcal A},
	}
	in which $P$ is as follows: 
	\begin{enumerate}
		\item When $n=2$ and ${\mathcal A}=\<\DS_2\>$, $P\in\left\{I_2, \mtx{1&0\\0 &-1} \right\}$. The two automorphisms are the identity map and the map
		\be{equation}{\label{<DS_2> nontrivial auto}
			\phi\left(\mtx{a &b\\b&a}\right)=\mtx{a &-b\\-b&a},\quad a, b\in\R.
		}
		\item When $n\ne 2$ or ${\mathcal A}\ne \<\DS_2\>$, $P\in {\mathcal A}$ is invertible. Moreover, when ${\mathcal A}=\<\DS_n\>$, 
		$P$ can be chosen from invertible matrices in $\DS_n$. 
		
	\end{enumerate}
\end{theorem}

\begin{proof} First, we prove the theorem for ${\mathcal A}=\<\DS_n\>$.
	\eqref{<DS> unitary} shows that the algebra $\<\DS_n\>\simeq \R\times M_{n-1}(\R)$  and
	an automorphism $\phi$ of $\<\DS_n\>$ has the form
	$$U\mtx{a&0\\0 &B}U^{-1}\overset{\phi}{\mapsto} U\psi\left(\mtx{a&0\\0 &B}\right) U^{-1}$$ in which
	$\psi$ is an automorphism of $\R\times M_{n-1}(\R)$. The automorphisms of $\R\times M_{n-1}(\R)$ must be direct sums of  irreducible representations of 
	$\R\times M_{n-1}(\R)$, which are $\R^1$ and $\R^{n-1}$ (cf. e.g. \cite[\S 3.3]{Etingof}).
	
	When $n=2$,
	the automorphism $\psi$  of the algebra $\R\times \R$ has two possibilities:    the identity map and the map   $\mtx{a&0\\0&b}\mapsto \mtx{b&0\\0&a}$.
	The corresponding automorphisms of $\<\DS_2\>$ are $\phi(A)=  PAP^{-1}$ for $P\in\left\{I_2, \mtx{1&0\\0 &-1} \right\}$. 
	
	When $n\ge 3$, by Skolem-Noether theorem (cf. \cite[Theorem 1]{Szi}), an algebra automorphism $\psi$ of $\R\times M_{n-1}(\R)$  has the following form for a given
	invertible $Q\in M_{n-1}(\R)$:
	\be{equation}{
		\psi\left(\mtx{a &0\\0 &B}\right)=\mtx{a &0\\0 &QBQ^{-1}},\qquad \mtx{a &0\\0 &B}\in \R\times M_{n-1}(\R).
	}
	The corresponding algebra automorphism $\phi$ of $\<\DS_n\>$ acts on $A=U\mtx{a&0\\0 &B}U^{-1}$ as
	\be{equation}{
		\phi(A)=\phi\left(U\mtx{a&0\\0 &B}U^{-1}\right)=U\mtx{a &0\\0 &QBQ^{-1}}U^{-1}.
	}
	Let $P_\epsilon=U\mtx{1&0\\0& \epsilon Q}U^{-1}$ for $\epsilon\in\R$. When $\epsilon\ne 0$, $P_\epsilon\in\<\DS_n\>$ is invertible, and 
	\be{equation}{
		\phi(A)=P_\epsilon AP_{\epsilon}^{-1}.
	}
	We have $\lim_{\epsilon\to 0} P_{\epsilon}=P_{0}=\frac{1}{n}ee^t$.
	Birkhoff's theorem (cf. \cite[Theorem 8.7.2]{Horn}) shows that $\DS_n$ is the convex hull of all $n\times n$ permutation matrices. The matrix $\frac{1}{n}ee^t$ is in the interior of
	$\DS_n$ in $\<\DS_n\>$. Hence $P_{\epsilon}\in\DS_n$ when $|\epsilon|\ne 0$ is sufficiently small.

	Second,  we prove   for ${\mathcal A}=\<\RS_n\>$. 
	By \eqref{<RS> unitary}, an automorphism $\phi$ of $\<\RS_n\>$ has the form
	\be{equation}{\label{RS automorphism 1}
		U\mtx{a&r^t\\0 &B}U^{-1}\overset{\phi}{\mapsto}
		U\psi\left( \mtx{a&r^t\\0 &B}\right) U^{-1}
	}
	in which $\psi$ is an automorphism of the algebra
	\be{equation}{
		{\mathcal A'}=\left\{\mtx{a&r^t\\0 &B} :\ a\in\R, B\in M_{n-1}(\R), r\in\R^{n-1}\right\}.
	}
	Skolem-Noether theorem implies that an automorphism of the algebra $M_{n-1}(\R)$ must 
	have the form $B\mapsto QBQ^{-1}$ for a fixed invertible $Q\in M_{n-1}(\R)$. 
	By \cite[Theorem 3.2]{Anh}, the automorphism $\psi$  is associate to  
	a vector  $v\in\R^{n-1}$ and invertible matrices $C, Q\in M_{n-1}(\R)$ such that
	\be{equation}{\label{RS block upper tri auto}
		\psi\left( \mtx{a&r^t\\0 &B}\right) = \mtx{a & av^t-v^tQBQ^{-1}+r^t C\\0 &QBQ^{-1}}
	}
	Applying \eqref{RS block upper tri auto} to the following identity for any $r\in\R^{n-1}$ and any invertible $B\in M_{n-1}(\R)$:
	\be{equation}{
		\psi\left( \mtx{1&r^t\\0 &B^{-1}}\right)\psi\left(\mtx{1&0\\0 &B}\right)=
		\psi\left( \mtx{1&r^t\\0 &B^{-1}}\mtx{1&0\\0 &B}\right) 
		=\psi\left( \mtx{1&r^tB\\0 &I_{n-1}}\right), 
	}
	we get $r^tCQBQ^{-1}=r^t BC$.  Since $r\in\R^{n-1}$ is arbitrary, $CQBQ^{-1}=BC$, that is, $BCQ=CQB$.
	Since $B$ is an arbitrary matrix, $CQ=cI_{n-1}$ for certain nonzero $c\in\R$.  \eqref{RS block upper tri auto} becomes
	\be{equation}{\label{RS block upper tri auto 2}
		\psi\left( \mtx{a&r^t\\0 &B}\right) = \mtx{a & av^t-v^tC^{-1}BC+r^t C\\0 &C^{-1}BC}
		={P'}^{-1} \mtx{a&r^t\\0 &B}{P'}
	}
	in which $P'=\mtx{1 &v^t\\0 &C}$. The automorphism $\phi$ of $\<\RS_n\>$ in \eqref{RS automorphism 1} becomes
	\be{equation}{
		A=U\mtx{a&r^t\\0 &B}U^{-1}\overset{\phi}{\mapsto}
		U {P'}^{-1} \mtx{a&r^t\\0 &B} P'U^{-1}
		= U {P'}^{-1} U^{-1} A U{P'} U^{-1}.
	}
	Let $P=U {P'}^{-1} U^{-1}\in \<\RS_n\>$. Then $\phi(A)=PAP^{-1}$ as desired. 
	
	The proof for ${\mathcal A}=\<\CS_n\>$ is similar. 
\end{proof}

We will need the following result regarding when a product  of stochastic matrices equals a permutation matrix.

\begin{lemma}\label{thm: product permutation}
	If $A_1,\ldots,A_m \in \RS_n$ (resp. $\DS_n$, $\CS_n$) and $A_1A_2\cdots A_m\in\P_n$, then each  $A_1,\ldots,A_m\in\P_n$.
\end{lemma}

\begin{proof}
	Without loss of generality, we  prove the statement for $A_1,\ldots,A_m \in \RS_n$. 
	The case of $\CS_n$ is analogous; and the case of $\DS_n$ can be derived from that of $\RS_n$.
	
	The statement obviously holds when $m=1$. 
	
	Now assume that $m>1$, let $A_1=[a_{ij}]\in\RS_n$ and $B=A_2\cdots A_m=[b_{ij}]\in \RS_n$. Each row sum of $A_1$  and row sum of $B$ is $1$.
	Suppose that 
	the $k$-th column of $B$ has the least column sum, which should be   no more than $1.$
	By assumption, $P=A_1B=[p_{ij}]$ is a permutation. Suppose $p_{\ell k}=1$ on the $k$-th column of $P$.
	Then 
	$$1=p_{\ell k}= \sum_{s=1}^{n} a_{\ell s}b_{sk} \le 
	(\sum_{s=1}^{n} a_{\ell s})(\sum_{s=1}^{n} b_{sk})= \sum_{s=1}^{n} b_{sk}\le 1.$$
	Therefore, the least column sum of $B$ is $\sum_{s=1}^{n} b_{sk}=1$, which implies that $B\in\DS_n$. Analogous argument shows that
	$A_1\in\DS_n$. 
	
	For each $i\in [n]$ there is a unique $j\in [n]$ such that 
	$$1=p_{ij}=\sum_{s=1}^{n} a_{is}b_{sj}\le (\sum_{s=1}^{n} a_{is})(\sum_{s=1}^{n} b_{sj})=1.$$
	The above equality holds if and only if there is $s\in [n]$ such that $a_{is}=b_{sj}=1,$ and
	$a_{ip}=b_{pj}=0$ for all $p\in[n]\setminus\{s\}$, if and only if both $A_1$ and $B=A_2\cdots A_m$ are permutation matrices. 
	Applying induction, we see that  $A_2,\ldots,A_m\in\P_n$. 
\end{proof}

\section{Multiplicative preservers on $\DS_n$ and $\<\DS_n\>$} 

In this section, we investigate the maps $\phi_i: \DS_n\to\DS_n$ ($i\in [m]$) and the maps $\phi_i: \<\DS_n\>\to\<\DS_n\>$ ($i\in [m]$) that preserve the multiplicative spectrum or multiplicative trace.
When $m=1$, it is too general to describe trace preservers, and  linearity should be assumed for spectrum preservers. 
When $m=2$, Theorem \ref{thm: two maps preserving trace} shows that multiplicative trace preservers $\phi_i$ can be extended to linear bijections on
$\<\DS_n\>$.  When $m\ge 3$, the  multiplicative spectrum  preservers and the  multiplicative  trace preservers turn out to be identical.

\subsection{Multiplicative spectrum and trace preservers on  $\DS_n$}

We  first describe the linear spectrum preservers on $\DS_n$ as follows. 

\begin{theorem}
	\label{thm: DS_n one map spectrum}
	A linear map   $\phi: \DS_n \to \DS_n$  preserves the spectrum:
	\begin{equation*}
		\spec(\phi(A))=\spec(A),\qquad A\in\DS_n,
	\end{equation*}
	if and only if there is a permutation $P\in\P_n$ such that
	\begin{equation}\label{DSn spec preserver}
		\phi(A)= PAP^t,\quad \text{or}  \quad \phi(A)= PA^t P^t. 
	\end{equation} 
\end{theorem}

\begin{proof} 
	Each permutation $A\in\DS_n$ satisfies $A^m=I_{n}$ for certain $m\in\Z_+$.
	Every linear spectrum preserver $\phi$ in $\DS_n$ should send the
	permutation $A\in\DS_n$ to $\phi(A)\in \DS_n$ such that $\spec(\phi(A)^m)=\spec(\phi(A))^m=\spec(A)^m$ has all eigenvalues 1, which means that $\phi(A)^m=I_{n}$.
	By Lemma \ref{thm: product permutation},
	$\phi(A)\in\P_n$. 
	
	If $A, B\in\P_n$ are distinct permutations such that $\phi(A)=\phi(B)$, then 
	$C=\frac{1}{2}(A+B)\in\DS_n$ is not a permutation but $\phi(C)=\phi(A)$ is a permutation. 
	There is $m'\in\Z_+$ such that $\phi(C)^{m'}=I_{n}$. Notice that 
	$$
	\spec(C^{m'})=\spec(C)^{m'}=\spec(\phi(C))^{m'} = \spec(\phi(C)^{m'}) =\spec(I_{n}).
	$$
	We have $C^{m'}=I_{n}$ so that $C$ is a permutation, which is a contradiction. 
	
	Now the linear map $\phi$ sends the set $\P_n$ bijectively to itself. 
	By \cite[Theorem 2.2]{LiTamTsing}, $\phi$ has one of the following forms for permutations $P,Q\in\P_n$:
	$$
	\phi(A)=PAQ,\qquad\text{or}\qquad \phi(A)=PA^t Q. 
	$$
	Using $\spec(\phi(I_{n}))=\spec(I_{n})$, we get \eqref{DSn spec preserver}. 
\end{proof}


\begin{theorem}
	\label{thm: DS_n two maps preserving trace of product}
	For two maps $\phi_i: \DS_n\to \DS_n$ ($i=1,2$), the following   are equivalent:
	\begin{enumerate}
		\item
		$\phi_1$ and $\phi_2$ satisfy that
		\be{equation}
		{\label{DS_n two map preserving product}
			\tr(\phi_1(A) \phi_2(B))=\tr(AB),\qquad  A, B\in \DS_n.
		}
		
		\item
		$\phi_1$ and $\phi_2$ satisfy that
		\be{equation}
		{\label{DS_n two map preserving product spectrum}
			\spec(\phi_1(A) \phi_2(B))=\spec(AB),\qquad  A, B\in \DS_n.
		}
		
		\item
		There exist permutation matrices $P, Q\in \P_n$ such that
		\be{subequations}{
			\be{eqnarray}
			{\label{DS_n two map forms 1}
				\phi_1(A)=PAQ^t,\quad \phi_2(A)=QAP^t,&& A\in\DS_n; \quad
				\text{or}
				\\ \label{DS_n two map forms 2}
				\phi_1(A)=PA^t Q^t,\quad \phi_2(A)=QA^t P^t,&& A\in\DS_n.
			}
		}
	\end{enumerate}
	If furthermore,  \eqref{DS_n two map preserving product} or \eqref{DS_n two map preserving product spectrum} holds and
	$\phi_1=\phi_2$, then $P=Q$ in \eqref{DS_n two map forms 1} and \eqref{DS_n two map forms 2}.
\end{theorem}

\begin{proof}
	The proofs of ``(3)$\Rightarrow$(2)$\Rightarrow$(1)'' are straightforward.
	We prove ``(1)$\Rightarrow$(3)'' here.
	
	Suppose $\phi_1$ and $\phi_2$ satisfy \eqref{DS_n two map preserving product}.
	By  Theorem \ref{thm: two maps preserving trace}, both $\phi_1$ and $\phi_2$ can be extended to linear bijections
	$\<\DS_n\>\to\<\DS_n\>$. 
	Given $A, B\in\DS_n$, we have $AB\in \DS_n$, so that
	$\tr(AB)\le n$; moreover, $\tr(AB)=n$ if and only if $A,B\in\P_n$ and $B=A^t$. 
	For any permutation $A\in\P_n$, \eqref{DS_n two map preserving product} implies that $\tr(\phi_1(A)\phi_2(A^t))=\tr(AA^t)=n$, so that 
	$\phi_1(A)\in\P_n$ and
	\be{equation}
	{\label{DS_n two map relation}
		\phi_2(A^t)=\phi_1(A)^t.
	}
	Therefore,  $\phi_1$ and $\phi_2$ are linear maps on $\DS_n$ that send  $\P_n$ to itself bijectively. By \cite[Theorem 2.2]{LiTamTsing},
	$\phi_1$ is of the form $\phi_1(X)=PXQ^t$ or $\phi_1(X)=PX^t Q^t$ for some  $P, Q\in\P_n$.  By \eqref{DS_n two map relation}
	and  Birkhoff's theorem, every doubly stochastic matrix is a convex combination of permutation matrices, we get
	$\phi_2(X)=\phi_1(X^t)^t$ for $X\in \DS_n$. Therefore, we get \eqref{DS_n two map forms 1} and \eqref{DS_n two map forms 2}. 
	
	The last claim for $\phi_1=\phi_2$ is obvious.
\end{proof}


\begin{theorem}
	\label{thm: DS_n m maps preserving trace of product}
	Let $m\ge 3$. 
	For  maps $\phi_i: \DS_n\to \DS_n$ ($i\in[m]$), the following   are equivalent:
	\begin{enumerate}
		\item
		$\phi_1,\ldots,\phi_m$ satisfy that
		\be{equation}
		{\label{DS_n m map preserving product}
			\tr(\phi_1(A_1) \phi_2(A_2)\cdots\phi_m(A_m))=\tr(A_1A_2\cdots A_m),\quad  A_1, \ldots, A_m\in \DS_n.
		}
		
		\item 
		$\phi_1,\ldots,\phi_m$ satisfy that
		\be{equation}
		{\label{DS_n m map preserving product spectrum}
			\spec(\phi_1(A_1) \phi_2(A_2)\cdots\phi_m(A_m))=\spec(A_1A_2\cdots A_m),\quad  A_1, \ldots, A_m\in \DS_n.
		}
		
		\item
		There exist permutation matrices $P_1, P_2,\ldots, P_m, P_{m+1}=P_1\in \P_n$ such that
		\be{equation}
		{\label{DS_n m maps preserve trace formula}
			\phi_i(A)=P_iAP_{i+1}^t,\quad  A\in\DS_n,\quad i\in [m].
		}
	\end{enumerate}
\end{theorem}

\begin{proof}
	The proofs of ``(3)$\Rightarrow$(2)$\Rightarrow$(1)'' are straightforward. 
	We prove ``(1)$\Rightarrow$(3)'' here.
	
	Suppose \eqref{DS_n m map preserving product} holds. 
	For $A_1, \ldots, A_m\in \DS_n$ we have $A_1\cdots A_m\in\DS_n$ so that
	$\tr(A_1\cdots A_m)\le n$. Moreover, by Lemma \ref{thm: product permutation}, $\tr(A_1\cdots A_m)= n$ if and only if $A_1,\ldots,A_m\in\P_n$   and 
	\be{equation}
	{\label{DS_n m maps trace max}
		A_1A_2\cdots A_m=A_2\cdots A_m A_1=\cdots=A_mA_1\cdots A_{m-1}=I_n.
	}
	For every $i\in [m]$ and $A_i\in\P_n$, there are permutation matrices $A_1,\ldots,A_{i-1}, A_{i+1},\ldots,A_m\in\P_n$ such that
	$A_1\cdots A_{i-1} A_iA_{i+1}\cdots A_m =I_n$. 
	By \eqref{DS_n m map preserving product},
	each $\phi_i$ ($i\in [m]$) maps $\P_n$ into $\P_n$. 
	
	Define $\wt{\phi_2}:\DS_n\to\DS_n$ such that
	\be{equation}
	{\label{DS_n wt phi_2}
		\wt{\phi_2}(B)=\phi_2(B)\phi_3(I_n)\cdots\phi_m(I_n).
	}
	Then  \eqref{DS_n m map preserving product} implies that
	\be{equation}
	{
		\tr(\phi_1(A)\wt{\phi_2}(B))=\tr(AB),\qquad  A, B\in \DS_n.
	}
	By Theorem \ref{thm: DS_n two maps preserving trace of product}, there exist  permutation matrices $P_1, P_2\in \P_n$ such that for $ A, B\in \DS_n$:
	\begin{subequations}
		\be{eqnarray}
		{\label{eq: DS_n m maps 1}
			\phi_1(A)=P_1AP_2^t,&& \wt{\phi_2}(B)=\phi_2(B)\phi_3(I_n)\cdots\phi_m(I_n)=P_2 B P_1^t;\qquad\text{or}
			\\ \label{eq: DS_n m maps 2}
			\phi_1(A)=P_1A^t P_2^t,&& \wt{\phi_2}(B)=\phi_2(B)\phi_3(I_n)\cdots\phi_m(I_n)=P_2 B^t P_1^t.
		}
	\end{subequations}
	
	First,  assume that \eqref{eq: DS_n m maps 1} holds.
	Since each   $\phi_3(I_n), \ldots, \phi_m(I_n) \in\P_n$, 
	$\phi_2(B)=P_2 B P_3^t$ for $B\in\DS_n$, where $P_3\in\P_n$. 
	For $i=2,\ldots,m$,  express \eqref{DS_n m map preserving product}  as
	\be{equation}
	{
		\tr(\phi_i(A_i) \cdots \phi_m(A_m) \phi_1(A_1)\cdots \phi_{i-1}(A_{i-1}))=\tr(A_i \cdots A_m A_1\cdots A_{i-1})
	}
	and apply the preceding argument. There exist permutation matrices $P_4,\ldots,P_m, P_{m+1}\in M_n$ such that
	$\phi_i(A_i)=P_i A_i P_{i+1}^t$ for $i=3,\ldots,m$.  \eqref{DS_n m map preserving product} implies that $P_{m+1}=P_1$.
	We get \eqref{DS_n m maps preserve trace formula}.
	
	Second, if \eqref{eq: DS_n m maps 2} holds, then similar argument shows that there exist  
	$P_3,\ldots,P_m\in \P_n$  such that
	$\phi_i(A_i)=P_i A_i^t P_{i+1}^t$ for $i=2,\ldots,m$.  Since $m\ge 3$, the trace
	$$
	\tr(\phi_1(A_1)\phi_2(A_2)\cdots\phi_m(A_m))=\tr(P_1 A_1^t A_2^t \cdots A_m^t P_{m+1}^t)
	$$
	is not always equal to $\tr(A_1\cdots A_m)$ for all $A_1,\ldots,A_m\in\DS_n$. This case is impossible.
\end{proof}

\subsection{Multiplicative spectrum and trace preservers on   $\<\DS_n\>$}

The linear spectrum preservers on $\<\DS_n\>$ can be described as follows. 

\begin{theorem}\label{thm: <DS> one map spectrum}
	A linear map   $\phi: \langle \DS_n\rangle \to \langle \DS_n\rangle $  preserves the spectrum:
	\begin{equation}
		\spec(\phi(A))=\spec(A),\qquad A\in \langle \DS_n\rangle,
	\end{equation}
	if and only if  $\phi$ has one of the following forms:
	\begin{enumerate}
		\item  $n=2$:   $\phi(A)=PAP^{-1}$ for  $P\in\left\{I_2, \mtx{1&0\\0 &-1} \right\}$. 
		\item $n\ge 3$: there is an invertible $P\in\<\DS_n\>$ such that 
		\begin{equation}\label{<DS_n> 1 map}
			\phi(A)=PAP^{-1} \qquad\text{or}\qquad 
			\phi(A)=PA^tP^{-1}.
		\end{equation}
	\end{enumerate}
\end{theorem}

\begin{proof}
	By Theorem \ref{thm: DS RS CS description}, 
	$\<\DS_n\>$ consists of $U\mtx{a&0\\0 &B}U^{-1}$ for $a\in\R$ and $ B\in M_{n-1}(\R)$. 
	A linear operator $\phi$ on $\langle \DS_n\rangle$ preserves the spectrum if and only if 
	the induced map $\wh{\phi}:\R\times M_{n-1}(\R)\to \R\times M_{n-1}(\R)$, defined by
	\begin{equation}\label{hat phi}
		\wh{\phi}\left(\mtx{a&0\\0 &B} \right) = U^{-1}\phi\left(U\mtx{a&0\\0 &B}U^{-1} \right)U,
	\end{equation}
	is a  linear spectrum preserver on   $\R \times M_{n-1}(\R)$. 
	It suffices to investigate the map $\wh{\phi}$. 
	\begin{enumerate}
		\item When $n=2$, the only spectrum preserver $\wh{\phi}$ sends $\diag(a,b)$ to either $\diag(a,b)$ or $\diag(b,a).$
		In consequence, we get $\phi(A)=A$ or $\phi(A)=\mtx{1 &0\\0 &-1} A \mtx{1 &0\\0 &-1} $ for $n=2$. 
		\item Suppose that $n\ge 3$.
		Write
		\begin{equation}
			\wh{\phi}(X)= \wh{\phi}_1(X)\oplus  \wh{\phi}_2(X)    
		\end{equation}
		for linear maps $\wh{\phi}_1:\R\times M_{n-1}(\R)\to \R$ and $\wh{\phi}_2:\R\times M_{n-1}(\R)\to M_{n-1}(\R).$ 
		
		Let $X=a\oplus B$ where $a\in\R$ and $B\in M_{n-1}(\R)$. On one hand, $\wh{\phi}_1(X)$ is in the spectrum of $\wh{\phi}(X)$, so that
		$\wh{\phi}_1(X)$ is in the spectrum of $X$. On the other hand,  $\wh{\phi}_1(X)$ is a linear combination of $a$ and the entries of $B\in M_{n-1}(\R)$.
		Since $n-1\ge 2$, the only possibility is $\wh{\phi}_1(X)=a$. 
		
		The operator
		$$\psi:M_{n-1}(\R)\to M_{n-1}(\R),\qquad
		\psi(A)=\wh{\phi}_2(0\oplus A).$$ is a linear spectrum preserver on $M_{n-1}(\R)$. Equivalently, the linear map $\psi$ satisfies the polynomial equation:
		\begin{equation}\label{spec poly preserver}
			\det(x I_{n-1}- A)=\det (x I_{n-1}-\psi(A)),\qquad A\in M_{n-1}(\R).
		\end{equation}
		We claim that the map 
		$\psi': M_{n-1}(\C)\to M_{n-1}(\C)$, defined by $$\psi'(A+\i B)=\psi(A)+\i\psi(B),\qquad A, B\in M_{n-1}(\R),$$ 
		is a linear spectrum preserver on $M_{n-1}(\C)$.  \eqref{spec poly preserver} implies that for any real number $y$ and any  $A, B\in M_{n-1}(\R)$, 
		$$\det(x I_{n-1}- A-yB)=\det (x I_{n-1}-\psi(A)-y\psi(B)).$$
		The two-variable polynomials on both sides of the above equality are identical on $(x,y)\in\C\times\R$. Thus they are identical on  $(x,y)\in\C\times\C$. In particular, 
		$$\det(x I_{n-1}- A-\i B)=\det (x I_{n-1}-\psi(A)-\i \psi(B)),$$
		which shows that $\psi'$  is a linear spectrum preserver on $M_{n-1}(\C)$.  By Marcus and Moyls' result in \cite{Ma59}, 
		either $\wh{\phi}_2\circ\iota(A)=SAS^{-1}$ or $\wh{\phi}_2\circ\iota(A)=SA^t S^{-1}$ for certain non-singular $S\in M_{n-1}(\R).$
		Consequently, the spectrum preserver $\wh{\phi}$ takes the form of either 
		\begin{equation}
			\wh{\phi}(a\oplus B) =a\oplus SBS^{-1} \qquad \text{or} \qquad  \wh{\phi}(a\oplus B) =a\oplus SB^t S^{-1} 
		\end{equation}
		for  certain non-singular $S\in M_{n-1}(\R).$ 
		We get   \eqref{<DS_n> 1 map} for $n\ge 3$. 
		\qedhere
	\end{enumerate}
	
\end{proof}

Theorem \ref{thm: two maps preserving trace} shows that every linear bijection  $\phi_1$ on   $\<\DS_n\>$ 
corresponds to a unique linear bijection $\phi_2$ on $\<\DS_n\>$ such that
$\tr(\phi_1(A_1)\phi_2(A_2))=\tr(A_1A_2)$ for $A_1, A_2\in\<\DS_n\>$; additional conditions should be added to obtain a refined characterization
that is equivalent to multiplicative spectrum preserver, as shown below.  

\begin{theorem}\label{thm: <DS> two maps spectrum trace}
	The following are equivalent for two maps $\phi_1, \phi_2: \<\DS_n\>\to \<\DS_n\>$:
	\begin{enumerate}
		
		\item 
		$\phi_1$ and $\phi_2$ satisfy that
		\be{equation}
		{\label{<DS_n> two map preserving trace}
			\tr(\phi_1(A_1) \phi_2(A_2))=\tr(A_1A_2),\quad  A_1, A_2\in \<\DS_n\>.
		}
		Moreover, there is   $p\in \{1,2\}$
		such that $A\in\<\DS_n\>$ and $\rank(A)\le 1$ imply that $\rank(\phi_p(A))\le 1$. 
		
		\item 
		$\phi_1$ and $\phi_2$ satisfy that
		\begin{equation}
			\label{<DS_n> two map preserving spectrum}
			\spec(\phi_1(A_1) \phi_2(A_2))=\spec(A_1A_2),\quad  A_1, A_2\in \<\DS_n\>.
		\end{equation}
		
		\item
		One of the following statements holds:
		\begin{enumerate}
			\item When $n=2$, there are invertible $C\in\<\DS_2\>$ and $P\in\left\{I_2, \mtx{1&0\\0 &-1} \right\}$ such that
			\be{equation}
			{\label{<DS_n> two map forms n=2} 
				\phi_1(A)=CPAP^{-1},\quad \phi_2(A)=PAP^{-1}C^{-1},\quad  A\in \<\DS_2\>.
			}
			
			\item When $n\ge 3$, there exist invertible matrices $P, Q\in \<\DS_n\>$ such that
			\begin{subequations}
				\be{eqnarray}
				{\label{<DS_n> two map forms n ge 3 a}
					\phi_1(A)=PAQ^{-1},\quad \phi_2(A)=QAP^{-1},&& A\in \<\DS_n\>;\quad  \text{or}
					\\ \label{<DS_n> two map forms n ge 3 b}
					\phi_1(A)=PA^t Q^{-1},\quad \phi_2(A)=QA^t P^{-1},&&   A\in \<\DS_n\>.
				}
			\end{subequations}
		\end{enumerate}
		
	\end{enumerate}
	
\end{theorem}

\begin{proof}
	The direction ``(3)$\Rightarrow$(2)''  is clear. 
	
	Next we prove ``(2)$\Rightarrow$(1)''. It is clear that \eqref{<DS_n> two map preserving spectrum} implies \eqref{<DS_n> two map preserving trace}. By Theorem \ref{thm: two maps preserving trace}, $\phi_1$ and $\phi_2$ are linear bijections. 
	Given $A\in \<\DS_n\>$ with $\rank(A)=1$, the bijectivity of $\phi_1$ and $\phi_2$ imply that $\phi_1(A)\ne 0$ and there exists $B\in\<\DS_n\>$ such that $\phi_2(B)=\phi_1(A)^t$. Then
	$$
	1\le \rank(\phi_1(A))=\rank(\phi_1(A)\phi_2(B))=\rank(AB)\le\rank(A)=1.  
	$$
	Therefore, $\rank(\phi_1(A))=1.$ The statement (1) is proved.

	Finally, we prove ``(1)$\Rightarrow$(3)''. Suppose \eqref{<DS_n> two map preserving trace} 
	holds, and  without lost of generality, $A\in\<\DS_n\>$ and $\rank(A)\le 1$  implies that $\rank(\phi_1(A))\le 1$.
	Let $U$ be given in \eqref{DS_n-orthogonal}.
	For $i=1,2$, define the maps
	$\wh{\phi}_i:\R\times M_{n-1}(\R)\to \R\times M_{n-1}(\R)$ such that
	\begin{equation}\label{hat phi 1}
		\wh{\phi}_i\left(\mtx{a&0\\0 &B} \right) = U^{-1}\phi_i\left(U\mtx{a&0\\0 &B}U^{-1} \right)U.
	\end{equation}
	Then  $\tr(\wh{\phi}_1(A_1)\wh{\phi}_2(A_2))=\tr(A_1A_2)$ and 
	$``\rank(A_1)\le 1 \implies \rank(\wh{\phi}_1(A_1))\le 1.''$  Both $\wh{\phi}_1$ and $\wh{\phi}_2$ are linear bijections
	by Theorem \ref{thm: two maps preserving trace}. We can retrieve $\phi_i$ from \eqref{hat phi 1} through the formula:
	\begin{equation}\label{nonhat phi 1}
		\phi_i(A)= U \wh{\phi}_i(U^{-1}AU) U^{-1},\quad A\in\<\DS_n\>. 
	\end{equation}
	\begin{enumerate}
		\item Suppose that $n=2$.  There are two possible forms for the bijective map $\wh{\phi}_1: \R\times\R\to \R\times \R$
		when $A\in\R\times\R$ and $\rank A\le 1$ implies that $\rank \wh{\phi}_1(A)\le 1$:  
		\begin{enumerate}
			\item $\wh{\phi}_1(\diag(a,b))=\diag(c_1a,c_2b)$ for some $c_1,c_2\in\R\setminus\{0\}$. 
			\item $\wh{\phi}_1(\diag(a,b))=\diag(c_1b,c_2a)$ for some $c_1,c_2\in\R\setminus\{0\}$.  
		\end{enumerate}  
		Using \eqref{nonhat phi 1}, we get two forms of $\phi_1$ described in \eqref{<DS_n> two map forms n=2} with the two possible $P$ matrices. When $\phi_1$ is fixed, \eqref{<DS_n> two map preserving trace} and Theorem \ref{thm: two maps preserving trace} imply  that
		$\phi_2$ is unique, which is  described in \eqref{<DS_n> two map forms n=2}. 
		
		\item Suppose that $n\ge 3$. We claim that 
		$$\wh{\phi}_1(\{0\}\times M_{n-1}(\R))\subseteq \{0\}\times M_{n-1}(\R).$$ 
		Suppose that on the contrary, there is $A'\in M_{n-1}(\R)$ such that $\wh{\phi}_1(0\oplus A')\not\in \{0\}\times M_{n-1}(\R)$. 
		Let $r=\rank(A')\ge 1$.   Then $A'=u_1v_1^t+\cdots+u_r v_r^t$ for some vectors $u_1,\ldots,u_r, v_1,\ldots,v_r\in\R^{n-1}$. 
		There is $i\in [r]$ such that  $\wh{\phi}_1(0\oplus u_iv_i^t)\not\in \{0\}\times M_{n-1}(\R)$. 
		By the rank constraint, $\wh{\phi}_1(0\oplus u_iv_i^t)$ must have rank 1.  
		So $\wh{\phi}_1(0\oplus u_iv_i^t)=a\oplus 0_{(n-1)\times (n-1)}$ for certain $a\ne 0$.
		Since $n\ge 3$, there exists $v_0\in\R^{n-1}$ that is not a scalar multiple of of $v_i$.
		The set $\{\wh{\phi}_1(0\oplus u_i(v_i+cv_0)^t)\mid c\in\R\}$ is a connected subset of $\R\times M_{n-1}(\R)$
		that contains $a\oplus 0_{(n-1)\times (n-1)}$ and consists of rank-1 matrices.
		Hence 
		$$\{\wh{\phi}_1(0\oplus u_i(v_i+cv_0)^t)\mid c\in\R\} \subseteq \R\times \{0_{(n-1)\times (n-1)}\}.$$
		The rank-1 matrices $\wh{\phi}_1(0\oplus u_iv_i^t)$ and $\wh{\phi}_1(0\oplus u_iv_0^t)$ are linear dependent.
		There is $c_0\in\R$ such that $\wh{\phi}_1(0\oplus u_i(v_i+c_0v_0)^t)=0$, which  contradicts  the bijectivity of 
		$\wh{\phi}_1$. So the claim is true. 
		
		As a linear bijection, $\wh{\phi}_1$ satisfies that  
		\begin{eqnarray*}
			\wh{\phi}_1(\{0\}\times M_{n-1}(\R)) &=& \{0\}\times M_{n-1}(\R), 
			\\ 
			\wh{\phi}_1(1\oplus 0_{(n-1)\times (n-1)}) &\not\in& \{0\}\times M_{n-1}(\R).
		\end{eqnarray*} 
		On one hand, since $\wh{\phi}_1$ preserves rank $\le 1$ matrices,
		$$\wh{\phi}_1(1\oplus 0_{(n-1)\times (n-1)})=c\oplus 0_{(n-1)\times (n-1)}$$ for some $c\ne 0$. 
		On the other hand, by Lim's description \cite{Lim} of linear preservers of  rank $\le 1$ matrices on $M_{m, n}(\F)$ for any field $\F$,   
		there are invertible matrices $M,N\in M_{n-1}(\R)$ such that 
		$$\wh{\phi}_1(0\oplus A)=0\oplus MAN
		\qquad\text{or}\qquad \wh{\phi}_1(0\oplus A)=0\oplus MA^t N.$$
		Consequently, using \eqref{nonhat phi 1},  $\phi_1$ has one of the forms described in
		\eqref{<DS_n> two map forms n ge 3 a} and \eqref{<DS_n> two map forms n ge 3 b}. 
		According to Theorem \ref{thm: two maps preserving trace}, when $\phi_1$ is fixed, the trace relation \eqref{<DS_n> two map preserving trace} uniquely determines   
		$\phi_2$, which also takes one of the forms in
		\eqref{<DS_n> two map forms n ge 3 a} and \eqref{<DS_n> two map forms n ge 3 b}.
		
	\end{enumerate}
	
	Overall, the proof is completed. 
\end{proof}

The multiplicative trace and spectrum  preservers  on $\<\DS_n\>$  for $m\ge 3$ maps are characterized as follows.

\begin{theorem}\label{thm: <DS> m ge 3 maps spectrum trace}
	Let $m\ge 3$. 
	The following   are equivalent for  maps $\phi_i: \<\DS_n\>\to \<\DS_n\>$ ($i\in [m]$):
	\begin{enumerate}
		\item
		$\phi_1,\ldots,\phi_m$ satisfy that
		\be{equation}
		{\label{<DS_n> m map preserving product}
			\tr(\phi_1(A_1) \phi_2(A_2)\cdots\phi_m(A_m))=\tr(A_1A_2\cdots A_m),\quad  A_1, \ldots, A_m\in \<\DS_n\>.
		}
		
		\item 
		$\phi_1,\ldots,\phi_m$ satisfy that
		\be{equation}
		{\label{<DS_n> m map preserving product spectrum}
			\spec(\phi_1(A_1) \phi_2(A_2)\cdots\phi_m(A_m))=\spec(A_1A_2\cdots A_m),\quad  A_1, \ldots, A_m\in \<\DS_n\>.
		}
		
		\item
		When $n=2$, there exist $P\in\left\{I_2, \mtx{1&0\\0 &-1} \right\}$ and invertible matrices $Q_1, Q_2,\ldots, Q_m\in\<\DS_2\>$ with $Q_1Q_2\cdots Q_m=I_2$ such that
		\be{equation}{\label{<DS_n> m maps preserve trace formula n=2}
			\phi_i(A)=PAP^{-1}Q_i,\quad  A\in \<\DS_2\>,\ i\in [m]. 
		}
		More explicitly, 
		\begin{subequations}
			\be{eqnarray}
			{\label{<DS_n> m maps preserve trace formula n=2 a}
				\phi_i(A)&=& AQ_i,\quad A\in \<\DS_2\>,\ i\in [m]; \quad \text{or}
				\\ \label{<DS_n> m maps preserve trace formula n=2 b}
				\phi_i\left(\mtx{a &b\\b&a}\right)&=&\mtx{a &-b\\-b&a}Q_i,\quad a, b\in\R,\ i\in [m].
			}
		\end{subequations}
		When $n\ge 3$, there are invertible matrices $P_1=P_{m+1}, P_2,\ldots, P_m $ in $\<\DS_n\> $ such that
		\be{equation}
		{\label{<DS_n> m maps preserve trace formula}
			\phi_i(A)=P_iAP_{i+1}^{-1},\quad  A\in \<\DS_n\>,\quad i\in [m].
		}
		
	\end{enumerate}
\end{theorem}

\begin{proof}
	The proofs of  ``(3)$\Rightarrow$(2)$\Rightarrow$(1)'' are straightforward.
	We now prove ``(1)$\Rightarrow$(3)''.  
	
	Let $\psi_2(B)=\phi_2(B)\phi_3(I_n)\cdots\phi_m(I_n)$. Then 
	\eqref{<DS_n> m map preserving product} implies that
	$\tr(\phi_1(A)\psi_2(B))=\tr(AB)$ for $A, B\in\<\DS_n\>$. 
	Theorem \ref{thm: two maps preserving trace} implies that $\phi_1$ and $\psi_2$ are linear bijections. So $\phi_i(I_n)$ are invertible for $i=3,\ldots,m$. 
	Similarly, each $\phi_i$ is a linear bijection  and each $\phi_i(I_n)$ is invertible  for $i\in [m]$. 
	
	Given $A_1, A_2\in\<\DS_n\>$, for all $A_3\in \<\DS_n\>$ we have
	\be{eqnarray}{
		&&\tr(\phi_3(A_3)\phi_4(I_n)\cdots\phi_m(I_n)\phi_1(A_1)\phi_2(A_2))
		\\
		&=&
		\tr(\phi_3(A_3)\phi_4(I_n)\cdots\phi_m(I_n)\phi_1(A_1A_2)\phi_2(I_n))
		\\ &=&
		\tr(\phi_3(A_3)\phi_4(I_n)\cdots\phi_m(I_n)\phi_1(I_n)\phi_2(A_1A_2)).
	}
	The surjectivity of $\phi_3$ and the invertibilities of $\phi_i(I_n)$ imply that 
	\be{equation}{\label{<DS> phi relations 1}
		\phi_1(A_1)\phi_2(A_2)=\phi_1(A_1A_2)\phi_2(I_n)=\phi_1(I_n)\phi_2(A_1A_2),
		\quad A_1, A_2 \in\<\DS_n\>. 
	}
	\eqref{<DS> phi relations 1} implies that:
	\begin{enumerate}
		\item [(I)]
		$\phi_1(I_n)\phi_2(A_2)=\phi_1(A_2)\phi_2(I_n)$, so that
		\be{equation}{\label{<DS_n> phi_2 phi_1}
			\phi_2(A_2)=\phi_1(I_n)^{-1}\phi_1(A_2)\phi_2(I_n),
			\quad
			A_2\in \<\DS_n\>.
		}
		
		\item [(II)]
		$\phi_1(A_1A_2)\phi_2(I_n)=\phi_1(A_1)\phi_2(A_2)=
		\phi_1(A_1)\phi_1(I_n)^{-1}\phi_1(A_2)\phi_2(I_n)$. Denote
		$\wt{\phi_1}(A)=\phi_1(A)\phi_1(I_n)^{-1}$. Then
		\be{equation}{
			\wt{\phi_1}(A_1A_2)=\wt{\phi_1}(A_1)\wt{\phi_1}(A_2),\quad
			A_1, A_2\in\<\DS_n\>.
		}
	\end{enumerate}
	
	Therefore, $\wt{\phi_1}$ is an algebra automorphism of $\<\DS_n\>$. 
	By Theorem \ref{thm: stochastic algebra automorphism},
	$\wt{\phi_1}(A)=PAP^{-1}$ in which $P$ is one of the following: 
	\begin{enumerate}
		\item [(I)]
		When $n=2$, we have  $P\in\left\{I_2, \mtx{1&0\\0 &-1} \right\}$ and
		\be{equation}{ \label{<DS_n> phi_1 expression}
			\phi_1(A)=\wt{\phi_1}(A)\phi_1(I_2)=PAP^{-1}\phi_1(I_2).
		}
		Matrices in $\<\DS_2\>$ are commutative. By \eqref{<DS_n> phi_2 phi_1} and \eqref{<DS_n> phi_1 expression},
		\be{eqnarray}{ \notag
			\phi_2(A) &=& \phi_1(I_2)^{-1}PAP^{-1}\phi_1(I_2) \phi_2(I_2)
			\\ \notag
			&=& P(P^{-1}\phi_1(I_2)^{-1}P)AP^{-1}\phi_1(I_2) \phi_2(I_2)
			\\ \notag
			&=& PA(P^{-1}\phi_1(I_2)^{-1}P)P^{-1}\phi_1(I_2) \phi_2(I_2)
			\\ \label{<DS_n> phi_2 expression}
			&=& PAP^{-1} \phi_2(I_2).
		}
		Likewise, $\phi_i (A) =PAP^{-1} \phi_i(I_2)$ for $i\in [m]$.  \eqref{<DS_n> m map preserving product} leads to
		$\phi_1(I_2)\cdots \phi_m (I_2)=I_2$. 
		We get \eqref{<DS_n> m maps preserve trace formula n=2},
		which is equivalent to \eqref{<DS_n> m maps preserve trace formula n=2 a} 
		and \eqref{<DS_n> m maps preserve trace formula n=2 b}.
		
		\item [(II)]
		When $n\ne 2$, $P\in \<\DS_n\>$ is invertible.  We have $\phi_1(A)=PA(\phi_1(I_n)^{-1}P)^{-1}$, and 
		\eqref{<DS_n> phi_2 phi_1} implies that
		$\phi_2(A)=(\phi_1(I_n)^{-1}P)A(\phi_2(I_n)^{-1}\phi_1(I_n)^{-1}P)^{-1}$. 
		We can  obtain identities  analogous to \eqref{<DS_n> phi_2 phi_1}:
		\be{equation}{\label{<DS_n> phi_(i+1) phi_i}
			\phi_{i+1}(A_{i+1})=\phi_i(I_n)^{-1}\phi_i(A_{i+1})\phi_{i+1}(I_n),
			\quad
			A_{i+1}\in \<\DS_n\>,\ i\in [m-1].
		}
		Let 
		\be{equation}{
			P_1=P; \quad P_{i+1}=\phi_{i}(I_n)^{-1}P_{i},\ i\in [m].
		}
		Then $\phi_{i}(A)=P_i A P_{i+1}^{-1}$ for $A\in\<\DS_n\>$ and $i\in [m]$. 
		\eqref{<DS_n> m map preserving product} implies that $P_{m+1}=P_1$. 
		\qedhere 
	\end{enumerate}
\end{proof}

The results in this section imply that every  multiplicative  spectrum (resp. trace) preserver $\phi_1,\ldots,\phi_m$  on  
$\DS_n$ can be extended  linearly  to a  multiplicative  spectrum (resp. trace) preserver on $\<\DS_n\>$, but the reverse is not necessarily true. 

\section{Multiplicative preservers on $\RS_n$ and $\<\RS_n\>$}


A matrix is in $\CS_n$  (resp.  $\<\CS_n\>$) if and only if its transpose is in $\RS_n$ (resp. $\<\RS_n\>$).  We will discuss
multiplicative  trace and spectrum  preservers on $\RS_n$ and $\<\RS_n\>$ in this section. Analogous results hold for 
multiplicative  trace and spectrum  preservers on $\CS_n$ and $\<\CS_n\>$.

\subsection{Multiplicative spectrum and trace preservers on   $\<\RS_n\>$}

Theorem \ref{thm: DS_n one map spectrum} implies the following result.

\begin{theorem}\label{thm: span RSn one map spectrum}
	A linear map   $\phi: \langle \RS_n\rangle \to \langle \RS_n\rangle $  preserves the spectrum:
	\begin{equation*}
		\spec(\phi(A))=\spec(A),\qquad A\in \langle \RS_n\rangle,
	\end{equation*}
	if and only if  $\phi$ has one of the following forms:
	\begin{enumerate}
		\item  $n=2$:   $\phi(A)_{\DS}=PA_{\DS}P^{-1}$ for  $P\in\left\{I_2, \mtx{1&0\\0 &-1} \right\}$. 
		\item $n\ge 3$: there is an invertible $P\in\<\DS_n\>$ such that 
		\begin{equation}\label{<RS_n> 1 map}
			\phi(A)_{\DS}=PA_{\DS}P^{-1} \qquad\text{or}\qquad 
			\phi(A)_{\DS}=P{A_{\DS}}^tP^{-1}.
		\end{equation}
	\end{enumerate}
\end{theorem}

\begin{proof}
	Suppose that $\phi: \langle \RS_n\rangle \to \langle \RS_n\rangle $  preserves the spectrum.
	Then for $A\in \<\RS_n\>$, 
	\begin{equation*}
		\spec(\phi(A)_{\DS} )=\spec(\phi(A))=\spec(A)=\spec(A_{\DS})  = \spec(\phi(A_{\DS}) ). 
	\end{equation*}
	The map $\psi:\<\DS_n\>\to \<\DS_n\>$ given by  $\psi(A)=\phi(A)_{\DS} $ is a linear spectrum preserver on $\<\DS_n\>$.  
	The possible forms of $\psi$ is described in Theorem \ref{thm: DS_n one map spectrum}. In particular, $\psi$ must be a linear bijection
	on $\<\DS_n\>$. We claim that   $\phi(A)_{\DS}= \phi(A_{\DS})_{\DS}$ for  all $A\in\<\RS_n\>$.  
	Suppose on the contrary, there is $A\in \<\RS_n\>$ such that $\phi(A-A_{\DS})_{\DS} =\phi(A)_{\DS}- \phi(A_{\DS})_{\DS}\ne 0$.
	Since $A-A_{\DS}\in R_n$ is nilpotent, $\phi(A-A_{\DS})_{\DS}$ is also nilpotent. There is a nilpotent $C\in\<\DS_n\>$ such that
	$C+\phi(A-A_{\DS})_{\DS}$ is not nilpotent. The bijectivity of $\psi$ implies that
	there is  $B\in\<\DS_n\>$ such that  $\psi(B)=\phi(B)_{\DS}=C$. $B$ is also nilpotent.   
	So $B+A-A_{\DS}$ is nilpotent; however, it has the same spectrum as $\phi(B+A-A_{\DS})_{\DS}= C+\phi(A-A_{\DS})_{\DS}$, which is not nilpotent.
	We get a contradiction. Hence $\phi(A)_{\DS}= \phi(A_{\DS})_{\DS}$  for  all $A\in\<\RS_n\>$. 
	The statements then follow from Theorem \ref{thm: DS_n one map spectrum}. 
\end{proof}

When $m=2$, Theorem \ref{thm: DS RS preserver relation} (2) and 
Theorem \ref{thm: <DS> two maps spectrum trace} directly lead to the following result.

\begin{theorem}\label{thm: <RS> two maps spectrum trace}
	The following statements are equivalent for two maps $\phi_i: \<\RS_n\>  \to \<\RS_n\>$ ($i=1,2$): 
	\begin{enumerate}
		\item 
		$\phi_1$ and $\phi_2$ satisfy that
		\be{equation}
		{\label{<RS_n> two map preserving spectrum}
			\spec(\phi_1(A_1) \phi_2(A_2))=\spec(A_1A_2),\quad  A_1, A_2\in \<\RS_n\>.
		}
		
		\item 
		$\phi_1$ and $\phi_2$ satisfy that
		\be{equation}
		{\label{<RS_n> two map preserving trace}
			\tr(\phi_1(A_1) \phi_2(A_2))=\tr(A_1A_2),\quad  A_1, A_2\in \<\RS_n\>.
		}
		Moreover, there is   $p\in \{1,2\}$ such that $A\in\<\DS_n\>$ and $\rank(A)\le 1$ imply that
		$\rank \phi_p(A)\le 1$. 
		
		\item 
		One of the following statements holds:
		\begin{enumerate}
			\item When $n=2$, there are invertible $C\in\<\DS_2\>$ and $P\in\left\{I_2, \mtx{1&0\\0 &-1} \right\}$ such that
			\be{equation}
			{\label{<RS_n> two map forms n=2} 
				\phi_1(A)_{\DS}=CPA_{\DS}P^{-1},\quad \phi_2(A)_{\DS}=PA_{\DS}P^{-1}C^{-1},\quad  A\in \<\RS_n\>.
			}
			
			\item When $n\ge 3$, there exist invertible matrices $P, Q\in \<\DS_n\>$ such that
			\begin{subequations}
				\be{eqnarray}
				{\label{<RS_n> two map forms n ge 3 a}
					\phi_1(A)_{\DS}=PA_{\DS}Q^{-1},\quad \phi_2(A)_{\DS}=QA_{\DS}P^{-1},&& A\in \<\RS_n\>;\quad  \text{or}
					\\ \label{<RS_n> two map forms n ge 3 b}
					\phi_1(A)_{\DS}=PA_{\DS}^t Q^{-1},\quad \phi_2(A)_{\DS}=QA_{\DS}^t P^{-1},&&   A\in \<\RS_n\>.
				}
			\end{subequations}
		\end{enumerate}
		
	\end{enumerate}
\end{theorem}

For $m\ge 3$, Theorem \ref{thm: DS RS preserver relation} (2) and 
Theorem \ref{thm: <DS> m ge 3 maps spectrum trace} imply the following result.

\begin{theorem}\label{thm: <RS> m ge 3 maps spectrum trace}
	When $m\ge 3$, 
	the following   are equivalent for  maps $\phi_i: \<\RS_n\>\to \<\RS_n\>$ ($i\in [m]$):
	\begin{enumerate}
		\item
		$\phi_1,\ldots,\phi_m$ satisfy that
		\be{equation}
		{\label{<RS_n> m map preserving product}
			\tr(\phi_1(A_1) \phi_2(A_2)\cdots\phi_m(A_m))=\tr(A_1A_2\cdots A_m),\quad  A_1, \ldots, A_m\in \<\RS_n\>.
		}
		
		\item 
		$\phi_1,\ldots,\phi_m$ satisfy that
		\be{equation}
		{\label{<RS_n> m map preserving product spectrum}
			\spec(\phi_1(A_1) \phi_2(A_2)\cdots\phi_m(A_m))=\spec(A_1A_2\cdots A_m),\quad  A_1, \ldots, A_m\in \<\RS_n\>.
		}
		
		\item
		When $n=2$, there exist $P\in\left\{I_2, \mtx{1&0\\0 &-1} \right\}$ and invertible matrices $Q_1, Q_2,\ldots, Q_m\in\<\DS_2\>$ with $Q_1Q_2\cdots Q_m=I_2$ such that
		\be{equation}{\label{<RS_n> m maps preserve trace formula n=2}
			\phi_i(A)_{\DS}=PA_{\DS} P^{-1}Q_i,\quad  A\in \<\RS_n\>,\ i\in [m]. 
		}
		When $n\ge 3$, there exist invertible matrices $P_1, P_2,\ldots, P_m, P_{m+1}=P_1\in \<\DS_n\> $ such that
		\be{equation}
		{\label{<RS_n> m maps preserve trace formula}
			\phi_i(A)_{\DS}=P_iA_{\DS} P_{i+1}^{-1},\quad  A\in \<\RS_n\>,\quad i\in [m].
		}

	\end{enumerate}
\end{theorem}

\subsection{Multiplicative spectrum and trace preservers on   $\RS_n$}

Theorem \ref{thm: span RSn one map spectrum} leads to the following result.

\begin{theorem}\label{thm: RS_n one map spectrum}
	A linear map   $\phi:  \RS_n  \to  \RS_n $  preserves the spectrum:
	\begin{equation*}
		\spec(\phi(A))=\spec(A),\qquad A\in  \RS_n ,
	\end{equation*}
	if and only if  $\phi$ has one of the following forms:
	\begin{enumerate}
		\item  $n=2$:   $\phi(A)= A_{\DS}+c A_{R}$ for a constant $c\in [-1, 1]$.
		\item $n\ge 3$: there is a permutation matrix $P\in \P_n$  such that
		\begin{equation}\label{RS one map spectrum preserver}
			\phi(A)=  PA  P^{-1} ,\qquad A\in\RS_n. 
		\end{equation}
		
	\end{enumerate}
\end{theorem}

\begin{proof} A map preserves the spectrum if and only if it preserves the characteristic polynomial. Since $\RS_n$ contains a nonempty open subset of $\<\RS_n\>$, every linear characteristic polynomial  preserver on   $\RS_n$ can be   extended to
	a linear characteristic polynomial  preserver on  $\<\RS_n\>$. Equivalently, every linear spectrum preserver on  $\RS_n$ can be uniquely extended to a linear spectrum preserver $\phi$  of $\<\RS_n\>$ such that  $\phi(\RS_n)\subseteq \RS_n.$      
	We discuss such $\phi$ case by case by applying  Theorem \ref{thm: span RSn one map spectrum}.  
	\begin{enumerate}
		\item Suppose that $n=2$. $\<\RS_2\>$ has a basis consisting of three matrices:
		$$A_0=\mtx{1 &-1\\1 &-1},\quad A_1=\mtx{1&0\\0&1},\quad \text{and} \ A_2=\mtx{0&1\\1&0}, $$
		in which $A_1$ and $A_2$ span $\<\DS_2\>$ and $A_0$ spans $R_2$. 
		The set $\RS_2$ is the convex combination of four matrices: $A_1$, $A_2$, and
		\begin{equation*}
			A_3=\mtx{1&0\\1&0},\quad A_4=\mtx{0&1\\0&1}.
		\end{equation*}
		We must have $\phi(A_i)\in\RS_2$ for $i=1,2,3,4$. The linear map $\wt{\phi}: \<\DS_2\>\to\<\DS_2\>$ defined by
		$\wt{\phi}(A)=\phi(A)_{\DS}$ is a spectrum preserver. By
		Theorem \ref{thm: span RSn one map spectrum}, there is  $P\in\left\{I_2, \mtx{1&0\\0 &-1} \right\}$ such that $\wt{\phi}(A)=PAP^{-1}$ for $A\in \<\DS_2\>$. 
		Since $\<\DS_2\>$ consists of real symmetric matrices, the only nilpotent matrix in $\<\DS_2\>$ is the zero matrix, which implies that  
		$\phi(A_0)_{\DS}=0$. Therefore, for $A\in\RS_n$, we have
		$\phi(A)_{\DS}=\phi(A_{\DS})_{\DS}=PA_{\DS}P^{-1}$, 
		so that $\phi(A)=PA_{\DS}P^{-1}+ g(A) A_0$ for a linear function $g:\<\RS_2\>\to\R$. 
		\begin{enumerate}
			\item     When $P=I_2$, we have $\phi(A)=A_{\DS}+g(A)A_0$. 
			The entries of $\phi(A_1)$ and $\phi(A_2)$ are nonnegative, so that
			$g(A_1)=g(A_2)=0$. In other words, $\<\DS_2\>\subseteq \ker g$.  Let $c=g(A_0)$. For $A\in\<\RS_2\>$, let $A_R=kA_0$ then 
			$$\phi(A)=A_{\DS}+g(A)A_0=A_{\DS}+g(A_R)A_0= A_{\DS}+ kg(A_0)A_{0}=A_{\DS}+cA_R. $$  
			Now $A_3=\frac{1}{2}\mtx{1&1\\1&1}+\frac{1}{2}\mtx{1&-1\\1&-1}$ and 
			$$\phi(A_3)= \frac{1}{2}\mtx{1&1\\1&1}+\frac{c}{2}\mtx{1&-1\\1&-1}\in\RS_2,$$
			which implies that $c\in [1,-1].$ The $c$ in this interval also makes $\phi(A_4)\in \RS_2.$
			
			\item When $P=\mtx{1&0\\0 &-1},$ we have $\phi(A_2)_{\DS}=P\mtx{0&1\\1&0}P^{-1}= \mtx{0&-1\\-1&0}.$ There is no matrix $\phi(A_2)_R\in R_2$ such that $\phi(A)=\phi(A_2)_{\DS}+ \phi(A_2)_R\in\RS_2$. This case is impossible. 
			
		\end{enumerate}
		
		\item Suppose that $n\ge 3$. Theorem \ref{thm: span RSn one map spectrum} shows that 
		$\phi(A)_{\DS}=PA_{\DS}P^{-1}$ or $\phi(A)_{\DS}=P{A_{\DS}}^tP^{-1}$ for an invertible $P\in\<\DS_n\>$.
		Consider $A\in\P_n$. There is $m$ such that $A^m=I$. So $(\phi(A)^{m})_{\DS}={\phi(A)_{\DS}}^{m}=PA^{m}P^{-1}=I$, which gives $\phi(A)^m=I$. 
		Since $\phi(A)\in\RS_n$, Lemma \ref{thm: product permutation} shows that $\phi(A)\in\P_n$. We get $\phi(\P_n)\subseteq \P_n.$
		
		\begin{enumerate}
			\item Suppose that  $\phi(A)_{\DS}=PA_{\DS}P^{-1}$ for an invertible $P\in\<\DS_n\>$.  Then
			$\phi(\P_n)= \P_n$.  By \cite[Theorem 2.2]{LiTamTsing}, we can choose $P\in\P_n$. 
			Consider linear spectrum preserver $\psi(A)=P^{-1}\phi(A)P$ on $\<\RS_n\>$, which   leaves $\RS_n$ invariant. 
			We have $\psi(A)=A$ for every $A\in \<\DS_n\>$. 
			
			We claim that $\psi(R_n)\subseteq R_n$. Otherwise, there exists a nonzero $B\in R_n$ such that $\psi(B)_{\DS}\ne 0$, which must be nilpotent as $B$ is nilpotent and $\spec(B)=\spec(\psi(B))=\spec(\psi(B))_{\DS}$. The nilpotent $\psi(B)_{\DS}$ cannot be skew-symmetric by spectral decomposition theorem. So the symmetric matrix ${\psi(B)_{\DS}}^t+\psi(B)_{\DS}\ne 0$, which is not nilpotent by spectral decomposition theorem.  However, ${\psi(B)_{\DS}}^t\in\<\DS_n\>$ so that
			$$
			{\psi(B)_{\DS}}^t+\psi(B)_{\DS}=\psi({\psi(B)_{\DS}}^t)+\psi(B)_{\DS}=\psi({\psi(B)_{\DS}}^t+B)_{\DS}.
			$$
			Since $\psi$ is a spectrum preserver, we get
			\begin{align*}
				\spec({\psi(B)_{\DS}}^t+\psi(B)_{\DS}) &=
				\spec (\psi({\psi(B)_{\DS}}^t+B)_{\DS})
				=
				\spec (\psi({\psi(B)_{\DS}}^t+B))
				\\
				&=\spec ({\psi(B)_{\DS}}^t+B)=\spec ({\psi(B)_{\DS}}^t)
				\\
				&=\spec ({\psi(B)_{\DS}})=\spec (\psi(B))=\spec(B),
			\end{align*}  
			which   contradicts the assumption that $B$ is nonzero nilpotent. Thus $\psi(R_n)\subseteq R_n$.
			
			Note that $R_n=\{ev^t\mid v\in\R^n,\ ve^t=0\}.$ If we write $\psi(ev^t)=ef(v^t)$ then $f$ is linear. There is $C=[c_{i,j}]_{n\times n}\in M_n$ such that $\psi(ev^t)=ev^t C$ for $ev^t\in R_n$.
			We get
			\begin{equation}
				\psi(A)=\psi(A_{\DS})+\psi(A_{R}) =  A_{\DS}  +A_{R}C.
			\end{equation}
			Consider $A_{1,2}=\mtx{1&0\\1&0}\oplus I_{n-2} \in \RS_n$. Then 
			$$(A_{1,2})_{\DS}= 
			\left (
			\begin{array}{cc|c}
				\frac{n-1}{n} & \frac{1}{n} & \multirow{2}{*}{ $O_{2 \times (n-2)}$} \\
				\frac{n-1}{n} & \frac{1}{n} & \\
				\hline
				-\frac{1}{n} & \frac{1}{n} & \\
				\vdots & \vdots & \multirow{2}{*}{ $I_{n-2}$} \\
				-\frac{1}{n} & \frac{1}{n} & 
			\end{array}
			\right),
			\quad
			(A_{1,2})_{R}=\mtx{1\\1\\\vdots\\1} \mtx{\frac{1}{n} &-\frac{1}{n} &0 &\cdots &0}. 
			$$
			The matrix $\psi(A_{1,2})= (A_{1,2})_{\DS} +(A_{1,2})_{R} C $ should be nonnegative. Checking the $(3,1)$    and the $(1,2)$  entries, we get
			\begin{equation}\label{A12 ineq}
				c_{1,1}-c_{2,1} \ge 1,\qquad
				c_{1,2}-c_{2,2} \ge -1.   
			\end{equation}
			Now consider $A_{2,1}=\mtx{0&1\\0&1}\oplus I_{n-2} \in \RS_n$. We have
			$$(A_{2,1})_{\DS}= 
			\left (
			\begin{array}{cc|c}
				\frac{1}{n} &\frac{n-1}{n} &  \multirow{2}{*}{ $O_{2 \times (n-2)}$} \\
				\frac{1}{n} &\frac{n-1}{n} &  \\
				\hline
				\frac{1}{n} &-\frac{1}{n} &  \\
				\vdots & \vdots & \multirow{2}{*}{ $I_{n-2}$} \\
				\frac{1}{n} & -\frac{1}{n} &  
			\end{array}
			\right),
			\quad
			(A_{2,1})_{R}=\mtx{1\\1\\\vdots\\1} \mtx{-\frac{1}{n} & \frac{1}{n} &0 &\cdots &0}. 
			$$
			Analogous analysis shows that
			\begin{equation}\label{A21 ineq}
				c_{2,1}-c_{1,1}\ge -1,\qquad 
				c_{2,2}-c_{1,2}\ge 1.   
			\end{equation}
			Comparing \eqref{A12 ineq} and \eqref{A21 ineq}, we get
			$c_{1,1}-c_{2,1}=1$ and $c_{2,2}-c_{1,2}=1$. 
			Similarly, we can get $c_{i,i}-c_{j,i}=1$ for  $i, j\in [n]$ and  $i\ne j$.
			Then $C=I_n + ew^t$ for certain $w\in\R^n$ so that for $A\in\<\DS_n\>$,
			$$
			\psi(A)=A_{\DS}+A_{R} (I_n + ew^t)=A_{\DS}+A_{R} =A.
			$$
			We get $\phi(A)=PAP^{-1}$ for $P\in\P_n$.

			\item Suppose that   $\phi(A)_{\DS}=P{A_{\DS}}^tP^{-1}$ for an invertible $P\in\<\DS_n\>$. Again  $\phi(\P_n)= \P_n$.   By \cite[Theorem 2.2]{LiTamTsing}, we get $P\in\P_n$. 
			Let  $\psi(A)=P^{-1}\phi(A)P$, which is a linear spectrum preserver   on $\<\RS_n\>$ that  leaves $\RS_n$ invariant. 
			Analogous to the preceding argument, we have $\psi(A)=A^t$ for $A\in\<\DS_n\>$ and $\psi(R_n)\subseteq R_n$, so that
			there is $C=[c_{i,j}]_{n\times n}\in M_n$ such that for $A=A_{\DS}+A_{R} \in\<\RS_n\>$, 
			\begin{equation}
				\psi(A)=\psi(A_{\DS})+\psi(A_{R}) =  {A_{\DS}}^t  +A_{R}C.
			\end{equation}
			On one hand, consider $A_{1,2}=\mtx{1&0\\1&0}\oplus I_{n-2} $  given in the preceding paragraph. 
			The matrix $\psi(A_{1,2})= {(A_{1,2})_{\DS}}^t  +(A_{1,2})_{R}C$ is nonnegative. 
			Checking the $(1,3)$ entry of $\psi(A_{1,2})$, we get $c_{1,3}-c_{2,3}\ge 1$. On the other hand, 
			for $A_{2,1}=\mtx{0&1\\0&1}\oplus I_{n-2}$, the image 
			$\psi(A_{2,1})=   {(A_{2,1})_{\DS}}^t  +(A_{2,1})_{R}C$ should be nonnegative. 
			The $(1,3)$ entry of $\psi(A_{2,1})$ shows that $c_{2,3}-c_{1,3}\ge 1$, which contradicts $c_{1,3}-c_{2,3}\ge 1$. 
			Therefore, this case is impossible. 
			\qedhere 
			
		\end{enumerate}
	\end{enumerate} 
\end{proof}

The multiplicative spectrum (resp. trace) preserver on   $\RS_n$ and $\CS_n$ are more sophisticated
than those of $\<\RS_n\>$ and $\<\CS_n\>$.

\begin{example} \label{RS_n: preserver in <RS_n> not RS_n}
	Let  $n\ge 3$. The maps $\phi_i(A) = A^t$ ($i=1,2$) form a multiplicative spectrum and   trace preserver on  $\DS_n$.
	We claim that $\phi_i$ ($i=1,2$) cannot be extended to a multiplicative spectrum and   trace preserver on  $\RS_n$. 
	
	Suppose on the contrary, these maps can be extended to 
	(with ambiguity of notions) $\phi_i: \RS_n\to\RS_n$ that  preserve   multiplicative spectrum and  trace. 
	By linearity of the trace,  $\phi_i$ can be extended to $\phi_i: \<\RS_n\>\to \<\RS_n\>$  ($i=1,2$) that preserve   multiplicative trace. 
	Moreover, since the trace operator is non-degenerate on $\DS_n$, Theorem \ref{thm: two maps preserving trace} shows that the extension from $\DS_n$ to $\<\DS_n\>$ is unique, which is $\phi_i(A)=A^t$ for $A\in \<\DS_n\>$.

	We have $\<\RS_n\>=\<\DS_n\>\oplus R_n$. 
	For every $A_1\in\<\DS_n\>$ and every $A\in \<\RS_n\>$  we have $\phi_1(A_1), \phi_2(A_{\DS})\in\<\DS_n\>$ and 
	\begin{align*}
		\tr(\phi_1(A_1)\phi_2(A_{\DS} )) &= \tr(A_1 A_{\DS}) =\tr(A_1A )
		\\ &
		= \tr(\phi_1(A_1)\phi_2(A ))
		=\tr(\phi_1(A_1)\phi_2(A )_{\DS}).
	\end{align*}
	Hence $\phi_2(A )_{\DS}=\phi_2(A_{\DS} ) ={A_{\DS}}^t$ and
	$$\phi_2(A )=\phi_2(A)_{\DS} + \phi_2(A)_R   =   {A_{\DS}}^t+\phi_2(A)_R  $$ 
	where $\gamma_2(A)=\phi_2(A)_R $ is a map from $\<\RS_n\>$ to $R_n$.
	Likewise, we get $\phi_1(A)={A_{\DS}}^t+\gamma_1(A )$ for a map $\gamma_1: \<\RS_n\> \to R_n.$
	
	Consider $B\in\RS_n$ such that
	\begin{align*}
		B &=\sum_{i=1}^{\lfloor\frac{n}{2}\rfloor}E_{i1}+\sum_{\lfloor\frac{n}{2}\rfloor+1}^{n} E_{j2}
		\ =\ B_{\DS}+B_{R},
		\\
		B_{R}&=  ev^t \in R_n,\qquad v^t=\big(\frac{\lfloor\frac{n}{2}\rfloor-1}{n}, \frac{\lfloor \frac{n+1}{2}\rfloor-1}{n}, -\frac{1}{n},\ldots, -\frac{1}{n}\big),
		\\
		B_{\DS} &=  B-ev^t \in \<\DS_n\>.
	\end{align*}
	For $j\in\{\lfloor\frac{n}{2}\rfloor+1,\ldots,n\}$, the $(j,1)$ entry of $B_{\DS}$ is $\frac{1-\lfloor\frac{n}{2}\rfloor}{n}$.
	For $i\in \{1,\ldots, \lfloor\frac{n}{2}\rfloor\}$, the $(i,2)$ entry of $B_{\DS}$ is $ \frac{1-\lfloor \frac{n+1}{2}\rfloor}{n}$. 
	We   claim that there is no map $\gamma:\RS_n\to R_n$ such that  ${B_{\DS}}^t+\gamma(B)\in\RS_n$. 
	Suppose on the contrary, such $\gamma$ exists. Then $\gamma(B)=e w^t$ in which $w^t=(w_1,\ldots,w_n)$ satisfies that
	$w^t e=0$. 
	Each entry of ${B_{\DS}}^t+e w^t$ should be nonnegative, which implies that
	\be{equation}{
		w_i\ge  \frac{ \lfloor \frac{n+1}{2}\rfloor-1}{n},\quad i\in \{1,\ldots, \lfloor\frac{n}{2}\rfloor\};
		\qquad 
		w_j\ge \frac{\lfloor\frac{n}{2}\rfloor-1}{n},\quad j\in\{\lfloor\frac{n}{2}\rfloor+1,\ldots,n\}.
	}
	Therefore, by $n\ge 3$,
	$$w^t e=\sum_{k=1}^{n} w_k\ge  \frac{ \lfloor \frac{n+1}{2}\rfloor-1}{n}\lfloor\frac{n}{2}\rfloor+\frac{\lfloor\frac{n}{2}\rfloor-1}{n} \lfloor \frac{n+1}{2}\rfloor
	=\be{cases}{m-1>0, &n=2m\\ m-2+\frac{m-1}{2m-1}>0, &n=2m-1}
	$$
	This contradicts the assumption $w^t e=0$. 
	Thus it is impossible to have  $\gamma:\RS_n\to R_n$ such that  ${B_{\DS}}^t+\gamma(B)\in\RS_n$. 
	This shows that $\phi_i(A)=A^t$ ($i=1,2$) cannot be extended to two maps on $\RS_n$ that form a trace or spectrum preserver.  
\end{example}

\begin{theorem}
	\label{thm: RS_n maps preserving trace of products}
	Let $m\ge 2$. The following are equivalent for maps  $\phi_i: \RS_n\to \RS_n$ ($i\in [m]$).
	\begin{enumerate}
		\item\label{RS_n trace statement}
		$\phi_i$ ($i\in [m]$) satisfy that
		\be{equation}
		{\label{RS_n maps preserving trace of product}
			\tr(\phi_1(A_1)\cdots\phi_m(A_m))=\tr(A_1\cdots A_m),\quad  A_1, \ldots, A_m\in \RS_n.
		}
		
		\item\label{RS_n spectrum statement}
		$\phi_i$ ($i\in [m]$) satisfy that
		\be{equation}
		{\label{RS_n maps preserving spectrum of product}
			\spec(\phi_1(A_1)\cdots\phi_m(A_m))=\spec(A_1\cdots A_m),\quad  A_1, \ldots, A_m\in \RS_n.
		}
		
		\item \label{RS_n trace spectrum formulas}
		There exist $n\times n$ permutation matrices $P_1 , P_2, \ldots, P_m, P_{m+1}=P_1$ such that for $ A=A_{\DS}+A_{R} \in\RS_n$: 
		\begin{equation}
			\label{RS_n m maps decomps}
			\phi_i(A)=P_i A_{\DS}P_{i+1} ^t+\gamma_i(A), \qquad i\in [m],
		\end{equation}
		in which each $\gamma_i: \RS_n\to R_n$ ($i\in [m]$) makes $\phi_i(A)\in\RS_n$ for all $A\in\RS_n$.
		Moreover, if $\phi_{\ell}$ is linear for certain $\ell\in [m]$, then 
		\begin{enumerate}
			\item when $n=2,$ there  is a  constant $c_{\ell}\in [-1,1]$ such that 
			\be{equation}
			{\label{RS_n m maps formula n = 2}
				\phi_{\ell}(A) =P_{\ell} (A_{\DS}+c_{\ell} A_{R})P_{\ell+1}^t=P_{\ell} A_{\DS}P_{\ell+1}^t+c_{\ell} A_{R}P_{\ell+1}^t, \quad  A\in \RS_2;
			}
			\item when $n\ge 3$,   
			\be{equation}
			{\label{RS_n m maps formula n ge 3}
				\phi_{\ell}(A) =P_{\ell} A P_{\ell+1}^t,\qquad  A\in \RS_n.
			}
		\end{enumerate}
	\end{enumerate}
\end{theorem}

\begin{proof} 
	The proof of ``\eqref{RS_n spectrum statement}$\Rightarrow$\eqref{RS_n trace statement}'' is easy. 
	
	``\eqref{RS_n trace spectrum formulas}$\Rightarrow$\eqref{RS_n spectrum statement}'': 
	By \eqref{DS parts preserve spectrum}, 
	for $A_1,\ldots,A_m\in\RS_n$:
	\be{equation}{
		\spec(\phi_1(A_1)\cdots\phi_m(A_m))=
		\spec(\phi_1(A_1)_{\DS}\cdots\phi_m(A_m)_{\DS})=
		\spec(A_1\cdots A_m).
	}
	Therefore, \eqref{RS_n spectrum statement} holds.

	``\eqref{RS_n trace statement}$\Rightarrow$\eqref{RS_n trace spectrum formulas}'': 
	Now suppose \eqref{RS_n maps preserving trace of product} holds. 
	Let $\RS_n^{(1)}$ be the subset of $\RS_n$ consisting of all matrices in which each row has exactly one nonzero entry 1. 
	Given $A, B\in\RS_n$, we have $AB\in \RS_n$, so that
	$\tr(AB)\le n$. 
	Moreover, $\tr(AB)= n$ if and only if $A, B\in \RS_n^{(1)}$ and $A=B^t$, if and only if $A=B^t$ is a permutation matrix.
	When $i\in[m]$ and $A_i\in\P_n$, for each $j\in [m]\setminus\{i\}$ we
	let $A_j={A_i}^t$ and $A_k=I_n$ for   $k\in [m]\setminus\{i,j\}$, then
	\eqref{RS_n maps preserving trace of product} and the preceding argument imply that  $\phi_i(A_i)\in\P_n$.

	Define
	$\phi'_i: \DS_n\to\<\DS_n \>$ ($i\in[m]$) such that
	\be{equation}
	{
		\phi'_i (A)=\phi_i(A)_{\DS},\quad   A\in\DS_n.
	}
	Theorem \ref{thm: DS RS CS} and \eqref{RS_n maps preserving trace of product} imply that 
	\be{equation}{\label{RS_n phi' trace eq}
		\tr(\phi_1'(A_1)\cdots\phi_m'(A_m))=\tr(A_1\cdots A_m),\quad A_1,\ldots,A_m\in \DS_n.
	}
	Set $\phi_{m+k}=\phi_k$ for $k\in [m-1]$. Define $\phi_j'':\DS_n\to\<\DS_n\>$ such that
	\be{equation}{
		\phi_j''(B)=(\phi_{j} (B)\phi_{j+1} (I_n)\cdots \phi_{j+m-2} (I_n))_{\DS}.
	}
	By Theorem \ref{thm: DS RS CS} and \eqref{RS_n maps preserving trace of product}, $\tr(\phi_i'(A) \phi_{i+1}''(B))=\tr(AB)$ for $A, B\in \DS_n$.  By Theorem \ref{thm: two maps preserving trace},
	$\phi_i'$  
	and $\phi_{i+1}''$ 
	can be extended to linear bijections $\psi_i'$ and $ \psi_{i+1}''$ on $\<\DS_n\>$, respectively, such that 
	\be{equation}{\label{RS_n psi' psi'' identity}
		\tr(\psi_i'(A) \psi_{i+1}''(B))=\tr(AB),\quad A, B\in \<\DS_n\>.
	}
	Moreover, for any $n\times n$ permutation matrix $P$, $\psi_i'(P)=\phi_i'(P)=\phi_i(P)_{\DS}$ is a permutation matrix. 
	Birkhoff's theorem shows that   $A\in\DS_n$ if and only if it is a convex combination of some permutation matrices $Q_1,\ldots,Q_N\in \P_n$:
	\be{equation}{
		A=t_1Q_1+\cdots+t_N Q_N,\qquad t_1,\ldots,t_N\in[0,1],\quad t_1+\cdots+t_N=1.
	} 
	Hence $\psi_i'(\DS_n) =\DS_n$. Likewise, $\psi_{i+1}''(\DS_n)=\DS_n$. By 
	\eqref{RS_n psi' psi'' identity}, Theorems \ref{thm: DS_n two maps preserving trace of product} 
	and   \ref{thm: DS_n m maps preserving trace of product},  there exist $n\times n$ permutation matrices $P_1,\ldots,P_m, P_{m+1}=P_1$ such that:
	\begin{enumerate}
		\item when $m=2$, 
		\be{eqnarray}
		{\label{RS_n two map forms}
			\text{either}&&
			\psi_1'(A)=P_1 A P_{2}^t,\quad \psi_2'(A)=P_2 A P_{1}^t,\quad   A\in \<\DS_n\>,
			\\\label{RS_n two map forms false}
			\text{or}&&
			\psi_1'(A)=P_1 A^t P_{2}^t,\quad \psi_2'(A)=P_2 A^t P_{1}^t,
			\quad   A\in \<\DS_n\>.
		}
		However, Example \ref{RS_n: preserver in <RS_n> not RS_n} shows that \eqref{RS_n two map forms false} is impossible for $n\ge 3$, and it is the same as \eqref{RS_n two map forms} when $n\le 2$.
		\item when $m\ge 3$, 
		\be{equation}
		{\label{RS_n m maps preserve trace formula}
			\psi_i'(A)=P_iAP_{i+1}^t,\quad  A\in\<\DS_n\>,\ i\in [m].
		}
	\end{enumerate}
	Moreover, by Theorem \ref{thm: DS RS CS} and \eqref{RS_n psi' psi'' identity}, for $A\in\RS_n$ and $B\in \DS_n$: 
	\be{eqnarray}{ \notag
		\tr(\phi_i(A)_{\DS}\psi_{i+1}''(B))
		&=&\tr(\phi_i(A)_{\DS}(\phi_{i+1}(B)\phi_{i+2}(I_n)\cdots \phi_{i+m-1}(I_n))_{\DS})
		\\\notag
		&=& \tr(\phi_i(A)\phi_{i+1}(B)\phi_{i+2}(I_n)\cdots \phi_{i+m-1}(I_n))
		\\\notag
		&=& \tr(AB)=\tr(A_{\DS}B)
		\\
		&=& \tr(\psi_i'(A_{\DS})\psi_{i+1}''(B)).
	}
	Since $\{\psi_{i+1}''(B): B\in\DS_n\}$ spans $\<\DS_n\>$, we have $\phi_i(A)_{\DS}=\psi_i'(A_{\DS})$ for $A\in\RS_n$. 
	Therefore, for $A\in\RS_n$, $\phi_i(A)$ has
	the possible form  \eqref{RS_n m maps decomps}, and
	$\gamma_i(A)=\phi_i(A)_{R}$ should let $\phi_i(A)\in\RS_n$.

	Now suppose furthermore  $\phi_{\ell}$ is linear for certain $\ell\in [m]$.
	Assume $\ell=1$ without loss of generality.  
	Then $\phi_1$ can be extended to a linear map $\wt{\phi_1}:\<\RS_n\>\to\<\RS_n\>$ such that
	\be{equation}{
		\wt{\phi_1}(A)=\wt{\phi_1}(A_{\DS})+\wt{\phi_1}(A_{R}),\quad A=A_{\DS}+A_{R}\in \<\RS_n\>.
	}  
	By preceding argument, $\phi_1$ permutes $n\times n$ permutation matrices. So $\phi_1$ maps $\DS_n$ bijectively and linearly onto $\DS_n$.
	When $A\in\DS_n$, 
	$\wt{\phi_1}(A_{\DS})=\phi_1(A)=\phi_1(A)_{\DS}=\psi_1'(A_{\DS})$. The linearity of $\wt{\phi_1}$ implies that 
	for all $A\in\RS_n$, $\wt{\phi_1}(A_{\DS})=\psi_1'(A_{\DS})$, which is described by \eqref{RS_n two map forms} and \eqref{RS_n m maps preserve trace formula}, and $\wt{\phi_1}(A_{R})=\phi_1(A)_R=\gamma_1(A)$. We conclude that both $\<\DS_n\>$ and $R_n$ are $\wt{\phi_1}$-invariant.
	
	Since
	$R_n=\{ev^t:\ v\in\R^n,\ v^t e=0\}$,
	there is a matrix $W \in M_{n}(\R)$ 
	such that 
	\be{equation}
	{\label{eq: RS_n two maps linear phi_i}
		\wt{\phi_1} (ev^t)=ev^t W,\qquad  ev^t\in R_n.
	}
	$\wt{\phi_1} (R_n)\subseteq R_n$ implies that $v^t We=0$ whenever $v^t e=0$. 
	So
	$We=c e$ for some $c\in\R$, or equivalently $W\in\<\RS_n \>$.

	By replacing $(\phi_i)_{i=1}^{m}$ with $(P_i^t \phi_i P_{i+1})_{i=1}^{m}$, 
	it suffices to prove \eqref{RS_n m maps formula n = 2} and \eqref{RS_n m maps formula n ge 3}   for the case $\ell=1$ and $P_1=P_2=I_n$.  Given $A=A_{\DS}+A_R \in\RS_n$ in which $A_R=ev^t$ with $v^t e = 0$, \eqref{eq: RS_n two maps linear phi_i} implies that
	\be{equation}
	{\label{RS_n m=2 phi_1}
		\phi_1(A)=A_{\DS}+ev^tW.
	}
	By Birkhoff's Theorem,  the convex hull of $\RS_n^{(1)}$ is $\RS_n$.  Therefore, the linear map $\phi_1$  maps
	$\RS_n$ into $\RS_n$ if and only if $\phi_1$ maps   $\RS_n^{(1)}$ into $\RS_n$.
	We discuss when $\phi_1$ maps   $\RS_n^{(1)}$ into $\RS_n$. 
	Assume that
	$W=[w_{ij}]\in \RS_n.$ 
	The case $n=1$ is trivial.
	
	\begin{enumerate}
		\item Suppose $n=2$. Then $W\in\<\RS_2 \>$ satisfies that $w_{11}+w_{12}=w_{21}+w_{22}$. The four matrices $A=A_{\DS}+A_{R}$ in $\RS_2^{(1)}$ always have $A_{\DS}=A_{\DS}^t$. So the two forms of $\phi_1$ in \eqref{RS_n m=2 phi_1} coincide.
		Let $B=\mtx{1&0\\1&0}=\mtx{1/2 &1/2\\1/2 &1/2}+\mtx{1\\1}\mtx{1/2 &-1/2}$.  
		Then
		\be{equation*}
		{
			\phi_1(B)=\mtx{1/2 &1/2\\1/2 &1/2}+\mtx{1\\1}\mtx{1/2 &-1/2}W=\mtx{\frac{1+w_{11}-w_{21}}{2} &\frac{1-w_{11}+w_{21}}{2}\\\frac{1+w_{11}-w_{21}}{2} &\frac{1-w_{11}+w_{21}}{2}}.
		}
		Let $c:=w_{11}-w_{21}$. Then $\phi_1(B)\in\RS_2$  implies that $|c|\le 1$. Moreover, when $|c|\le 1$, for all $A=A_{\DS}+ev^t\in\RS_2$:
		$$
		\phi_1(A) =A_{\DS}+ev^t\mtx{c+w_{21} &w_{12}\\ w_{21} &w_{12}+c}=A_{\DS}+cev^t.
		$$
		By computation, the above $\phi_1$ always maps $\RS_2^{(1)}$ into $\RS_2$.
		So \eqref{RS_n m maps formula n = 2}   is proved.

		\item Suppose $n\ge 3$.
		Consider
		\be{equation}
		{\label{eq: RS_n two maps n ge 3 B}
			B=\mtx{1 &0 &0 &\cdots &0\\0 &1 &0 &\cdots &0\\\vdots &\vdots &\vdots &\ddots &\vdots\\0 &1 &0 &\cdots &0}
			=B_{\DS}+e v^t
			\in\RS_n^{(1)}
		}
		in which
		\be{equation}{
			B_{\DS}:=\mtx{1 &-\frac{n-2}{n}  &\frac{1}{n} &\cdots &\frac{1}{n}\\ 0 &\frac{2}{n} &\frac{1}{n} &\cdots &\frac{1}{n}\\ \vdots &\vdots &\vdots &\ddots &\vdots\\ 0 &\frac{2}{n}  &\frac{1}{n} &\cdots &\frac{1}{n}},\quad
			v^t := \mtx{0, &\frac{n-2}{n},  -\frac{1}{n}, \cdots, -\frac{1}{n}}.
		}
		For each $n\times n$ permutation matrix $Q'$, $BQ' =B_{\DS}Q' +ev^tQ'\in\RS_n^{(1)}$.
		We have
		\be{equation}
		{\label{RS_n two maps permutation fix 1}
			\phi_1(BQ')=B_{\DS}Q'+ev^tQ'W \in\RS_n.
		}
		Denote the permutation matrix 
		\be{equation}
		{\label{RS_n permutation Q}
			Q:=E_{11}+\sum_{i=2}^{n-1} E_{i,i+1}+E_{n2}.
		}
		The $(2,1)$ entry of each $\phi_1(BQ^{r})$ for $r=0, 1,\ldots,n-2$ must be nonnegative. 
		Through \eqref{RS_n two maps permutation fix 1}, we get the system of inequalities
		\be{equation}
		{\label{eq: RS_n two maps ineqs 1}
			\be{cases}
			{
				\frac{n-2}{n}w_{21}-\frac{1}{n}w_{31}-\cdots- \frac{1}{n}w_{(n-1)1}- \frac{1}{n}w_{n1} &\ge 0
				\\
				-\frac{1}{n}w_{21}+\frac{n-2}{n}w_{31}-\cdots- \frac{1}{n}w_{(n-1)1}- \frac{1}{n}w_{n1} &\ge 0
				\\
				\qquad \qquad \qquad \qquad \vdots  &\vdots
				\\
				-\frac{1}{n}w_{21}-\frac{1}{n}w_{31}-\cdots- \frac{1}{n}w_{(n-1)1}+\frac{n-2}{n}w_{n1} &\ge 0
			}
		}
		Adding the above inequalities together yields $0\ge 0$. So all equalities hold  in  \eqref{eq: RS_n two maps ineqs 1}. We get $w_{21}=w_{31}=\cdots=w_{n1}$.
		Similar arguments show that $w_{ij}=w_{kj}$ for any distinct numbers $i,j,k\in\{1,\ldots,n\}$. Since $W\in\<\RS_n \>$, we get $W=c  I_n+ew^t$ in which $c\in\R$ and  
		$w\in\R^n$. Hence  for   $A=A_{\DS}+ev^t\in\RS_n$ (where $v^te=0$):  
		$$\phi_1(A)=A_{\DS}+ev^t (cI_n+ew^t)=A_{\DS}+cev^t.$$
		For the matrix $B$ in \eqref{eq: RS_n two maps n ge 3 B}, the $(1,2)$ and $(1,3)$ entries of $\phi_1(B)$ must be nonnegative, which lead to the inequalities:
		$$
		-\frac{n-2}{n}+\frac{n-2}{n} c \ge 0,\qquad \frac{1}{n}-\frac{1}{n}c \ge 0.
		$$
		Thus $c=1$ and
		\be{equation}
		{
			\phi_1(A)=A_{\DS}+ev^t=A,\qquad  A\in\RS_n.
		}
		Hence   \eqref{RS_n m maps formula n ge 3} holds.
		
	\end{enumerate}
	Overall, the proof of ``\eqref{RS_n trace statement}$\Rightarrow$\eqref{RS_n trace spectrum formulas}''  is completed.
\end{proof}

\section*{Acknowledgment}
M.C. Tsai is supported by Taiwan NSTC grant 114-2115-M-027-001.

\end{document}